\pdfoutput=1
\RequirePackage{ifpdf}
\ifpdf 
\documentclass[pdftex]{sigma}
\else
\documentclass{sigma}
\fi

\usepackage[all]{xy}

\newcommand{\R}{\mathbb{R}}

\newcommand{\C}{\mathcal{C}}
\newcommand{\D}{\mathcal{D}}
\newcommand{\M}{\mathcal{V}}
\newcommand{\MM}{\mathcal{V}_{\textrm{\normalfont I}}}
\newcommand{\MMM}{\mathcal{V}_{\textrm{\normalfont II}}}
\newcommand{\CC}{\mathcal{C}^{1}}
\newcommand{\A}{\boldsymbol{A}}
\newcommand{\E}{\mathcal{E}}
\newcommand{\B}{\boldsymbol{B}}

\newcommand{\V}{{\mathcal{V}}}
\newcommand{\Ann}{\operatorname{Ann}}
\newcommand{\Gr}{\operatorname{Gr}}
\newcommand{\Fl}{\operatorname{Fl}}
\newcommand{\Span}[1]{\left\langle #1 \right\rangle}
\renewcommand{\P}{\mathbb{P}}
\newcommand{\id}{\mathrm{id}}
\newcommand{\rank}{\operatorname{rank}}
\newcommand{\Char}[1]{\operatorname{char}^\R_{#1}}
\newcommand{\Th}{^\textrm{th}}
\newcommand{\St}{^\textrm{st}}
\newcommand{\Rd}{^\textrm{rd}}
\newcommand{\Nd}{^\textrm{nd}}
\newcommand{\Ddue}{{\check{\mathcal{D}}}}
\newcommand{\w}{\nu}
\newcommand{\elle}{H}
\newcommand{\p}{\varrho}

\numberwithin{equation}{section}

\newtheorem{Theorem}{Theorem}[section]
\newtheorem{Corollary}[Theorem]{Corollary}
\newtheorem{Lemma}[Theorem]{Lemma}
\newtheorem{Proposition}[Theorem]{Proposition}
 { \theoremstyle{definition}
\newtheorem{Definition}[Theorem]{Definition}
\newtheorem{Example}[Theorem]{Example}
\newtheorem{Remark}[Theorem]{Remark} }

\begin{document}


\renewcommand{\thefootnote}{$\star$}

\newcommand{\arXivNumber}{1403.3521}

\renewcommand{\PaperNumber}{032}

\FirstPageHeading

\ShortArticleName{$3^{\rm rd}$ Order Monge--Amp\`ere Equations}

\ArticleName{Meta-Symplectic Geometry of $\boldsymbol{3^{\rm rd}}$ Order\\
 Monge--Amp\`ere Equations and their Characteristics\footnote{This paper is a~contribution to the Special Issue
on Analytical Mechanics and Dif\/ferential Geometry in honour of Sergio Benenti.
The full collection is available at \href{http://www.emis.de/journals/SIGMA/Benenti.html}{http://www.emis.de/journals/SIGMA/Benenti.html}}}

\Author{Gianni MANNO~$^\dag$ and Giovanni MORENO~$^\ddag$}

\AuthorNameForHeading{G.~Manno and G.~Moreno}

\Address{$^\dag$~Dipartimento di Scienze Matematiche ``G.L.~Lagrange'', Politecnico di Torino,\\
\hphantom{$^\dag$}~Corso Duca degli Abruzzi 24, 10129 Torino, Italy}
\EmailD{\href{mailto:giovanni.manno@polito.it}{giovanni.manno@polito.it}}

\Address{$^\ddag$~Institute of Mathematics, Polish Academy of Sciences,\\
\hphantom{$^\ddag$}~ul.~\'Sniadeckich~8, 00-656 Warsaw, Poland}
\EmailD{\href{gmoreno@impan.pl}{gmoreno@impan.pl}}
\URLaddressD{\url{https://www.impan.pl/en/sites/gmoreno/home}}

\ArticleDates{Received October 29, 2015, in f\/inal form March 16, 2016; Published online March 26, 2016}

\Abstract{This paper is a natural companion of [Alekseevsky D.V., Alonso Blanco R., Manno G., Pugliese F., \textit{Ann. Inst. Fourier (Grenoble)} \textbf{62}
 (2012), 497--524, arXiv:1003.5177],
generalising its perspectives and results to the context of third-order~(2D) {M}onge--{A}mp\`ere equations, by using the so-called ``meta-symplectic structure'' associated with the 8D prolongation~$M^{(1)}$ of a~5D contact manifold~$M$.
We write down a geometric def\/inition of a~third-order {M}onge--{A}mp\`ere equation in terms of a~(class of) dif\/ferential two-form on~$M^{(1)}$. In particular, the equations corresponding to decomposable forms admit a simple description in terms of
 certain three-dimensional distributions, which are made from the characteristics of the original equations.
We conclude the paper with a study of the intermediate integrals of these special {M}onge--{A}mp\`ere equations, herewith called of Goursat type.}


\Classification{53D10; 35A30; 58A30; 14M15}

\renewcommand{\thefootnote}{\arabic{footnote}}
\setcounter{footnote}{0}

\section{Introduction}

Classical Monge--Amp\`ere equations (MAEs with one unknown function and two independent variables) constitute a distinguished class of scalar
$2\Nd$ order (non-linear) PDEs owing to a~remarkable property: the totality of their characteristics, herewith called
\emph{characteristic cone}\footnote{The notion behind it is rather old, and may have appeared in other
guises someplace else.}, degenerates into the union of zero, one or two 2D planes (according to the elliptic, parabolic or hyperbolic character of the equation).

The primary motivation of this paper was to determine whether, and to what extent, such a~phenomenon occurs also in the context of $3\Rd$ order PDEs, where there can be found physically interesting analogues of the classical MAEs \cite{MR1925763,MR2369200,MR1397274,MR1958112,MR1372461}. To the authors' best knowledge, such PDEs were def\/ined by Boillat by using Lax's complete exceptionality \cite{MR1139843,MR1194520, MR0067317}.
We introduce now, in a friendly coordinate way, the notion of characteristic cone, whose intrinsic def\/inition will be given later on (see~\eqref{eq.char.cone.1} and~\eqref{eq.char.cone.2}).

\subsection*{The main idea: studying PDEs using the geometry of their characteristics}\label{secIntroCarCone}
Let
\begin{gather}\label{equation}
\E\colon \ F\big(x^1,x^2,u,p_1,p_2,\dots,p_{i_1\cdots i_l},\dots\big)=0 , \qquad l \leq k,
\end{gather}
be a scalar
$k\Th$ order PDE with one unknown function $u=u(x^1,x^2)$ and two independent variables $(x^1,x^2)$.
As usual (see \cite{MR1504329} and, more recently, e.g.,~\cite{MR2813504,MR2389645} and references therein)
in the geometric theory of PDEs, the variables $p_{i_1\cdots i_l}$ correspond to the partial derivatives $\frac{\partial^{l}u}{\partial x^{i_1}\cdots\partial x^{i_l}}$, with the indices $ i_1\leq i_2\leq \cdots\leq i_l$ ranging in $ \{1, 2\}$.

A \emph{Cauchy problem} is obtained by complementing \eqref{equation} with some
 initial conditions
\begin{gather}\label{Cauchy_problem}
f \big(X^{1}(t),X^{2}(t)\big) =U(t),\qquad
\frac{\partial^\ell f}{\partial z^{\ell}} \big(X^{1}(t),X^{2}(t)\big) =Q_{\ell}(t), \qquad \ell\leq k-1,
\end{gather}
where
\begin{gather}\label{Cauchy_data}
 X_1(t), \ X_2(t),\ U(t),\ Q_\ell(t)
\end{gather}
 are given functions, and $\frac{\partial}{\partial z}$ is the derivative along the \emph{normal direction}\footnote{We adopted the same notation used in~\cite{LucaChar}, where the reader may also f\/ind a gentle introduction to the theory of characteristics and singularities of non-linear PDEs.} of the curve $(X^1(t)$, $X^2(t))$. A solution of the Cauchy problem~\eqref{equation},~\eqref{Cauchy_problem} is a function $u=f(x^1,x^2)$ which, together with its derivatives, satisf\/ies both~\eqref{equation} and~\eqref{Cauchy_problem}.
Initial data~\eqref{Cauchy_data} can be used to construct the curve
\begin{gather}\label{Cauchy_datum}
\Phi(t)=\big(X^{1}(t),X^{2}(t),U(t),P_{1}(t),P_{2}(t),\dots,P_{i_1\cdots i_{l}}(t),\dots\big),\qquad l\leq k-1,
\end{gather}
in the space with coordinates $(x^1,x^2,u,\dots,p_{i_1\cdots i_l},\dots)$, $l\leq k-1$, which we can then interpret as a \emph{Cauchy datum} (see, e.g., \cite{MorenoCauchy} for a jet-theoretic treatment of the space of Cauchy data). If such curve is \emph{non-characteristic} for equation~\eqref{equation}, then, assuming~$F$ real analytic, Cauchy problem~\eqref{equation},~\eqref{Cauchy_problem} admits a~(locally) unique analytic solution. Hence, \emph{characteristic curves} play a crucial role in the analysis of Cauchy problems. They can be def\/ined as follows.

A curve \eqref{Cauchy_datum} is \emph{characteristic} for $\E$ at the point\footnote{The choice of notation ``$\overline{m}^{k-2}$'' and ``$\overline{m}^{k-1}$'' will be motivated later on.} $\overline{m}^{k-1}\in\E$, i.e., a point
\begin{gather}\label{punto}
\overline{m}^{k-1}=\big(\overline{x}^{1},\overline{x}^{2},\overline{u},\dots ,\overline{p}_{i_1\cdots i_{l}},\dots\big),\qquad l\leq k,
 \end{gather}
whose coordinates satisfy \eqref{equation},
 if the tangent vector $\w=\dot{\Phi}(t_0)$ at the point
\begin{gather}\label{eq.point.m.k.meno.2}
\overline{m}^{k-2}=\Phi(t_0)=\big(\overline{x}^{1},\overline{x}^{2},\overline
{u},\dots,\overline{p}_{i_1\cdots i_{l}},\dots\big),\qquad l\leq k-1,
\end{gather}
 such that
\begin{gather}\label{characteristic}
\sum_{\ell_1+\ell_2=k} (-1)^{\ell_1} \frac{\partial F}{\partial p_{{\underbrace{\mbox{\scriptsize $1\cdots 1$}}_{\ell_1}}{\underbrace{\mbox{\scriptsize $2\cdots 2$}}_{\ell_2}}}} \Bigg|_{\overline{m}^{k-1}}\big(\w^2\big)^{\ell_1}\big(\w^1\big)^{\ell_2}= 0,
\end{gather}
where
\begin{gather}\label{eq.char.dir}
\w=\w^1\bigg({\partial}_{x^{1}}+\overline{p}_{1}\partial_{u}+\hspace{-7mm}\sum_{\underset{ h\leq k}{\scriptscriptstyle 1\leq i_2\leq\dots\leq i_h\leq 2}}\hspace{-5mm}\overline
{p}_{1i_2\cdots i_h}\partial_{p_{i_2\cdots i_{h}}}\bigg)
+
\w^2\bigg({\partial}_{x^{2}}+\overline{p}_{2}\partial_{u}+\hspace{-7mm}\sum_{\underset{h\leq k}{\scriptscriptstyle 1\leq i_2\leq\dots\leq i_h\leq 2}}\hspace{-5mm}\overline
{p}_{2i_2\cdots i_h}\partial_{p_{i_2\cdots i_{h}}}\bigg).
\end{gather}
From \eqref{characteristic}
one can associate with any point \eqref{punto} of $\E$ a number $\leq k$ of
directions \eqref{eq.char.dir} (the polynomial \eqref{characteristic} might possess multiple and/or imaginary roots) in the space with coordinates $(x^1,x^2,u,\dots,p_{i_1\cdots i_l},\dots)$, $l\leq k-1$. So, if we keep the point~\eqref{eq.point.m.k.meno.2} f\/ixed and let the point \eqref{punto} vary in $\E$, the aforementioned directions form, in general, a~cone.

The set $ \M^\E_{\overline{m}^{k-2}}:=\{\textrm{directions $\w$ as in \eqref{eq.char.dir}}\,|\,\exists\, \overline{m}^{k-1}\in\E \textrm{ such that $(\w^1,\w^2)$ satisf\/ies \eqref{characteristic}}\}$ is what we call the \emph{characteristic cone} of $\E$ at $\overline{m}^{k-2}$, and
the union
\begin{gather}\label{eq.char.cone.1}
\M^\E:=\bigcup_{\overline{m}^{k-2}} \M^\E_{\overline{m}^{k-2}}
\end{gather}
is a geometric object naturally associated to $\E$ (see Section \ref{sezionedovesidefinisceilconocaratteristico} later on).

In this paper, we basically
study those PDEs that correspond to 3D distributions on the 8D f\/irst prolongation of a 5D contact manifold.
Since the correspondence equation$\leftrightarrow$distribution is essentially the same as in the case of $2\Nd$ order MAEs, the PDEs we are interested in are precisely the $3\Rd$ order analogues
 of the classical (i.e., second-order) bidimensional MAEs.

Our main motivation was the apparent lack of a systematic geometric analysis of $3\Rd$ order MAEs,
carried out in the same spirit as for classical MAEs. Much as in the geometric approach to classical MAEs it is convenient to exploit the contact/symplectic geometry underlying $2\Nd$ order PDEs \cite{MR2383541, MR2503974,MR2805306,MR2352610,MR781584,MR1458653}, to deal with $3\Rd$ order MAEs we shall make use of the prolongation of a contact manifold, equipped with its Levi form (see, e.g., \cite[Section~2]{MR2279499}), a structure known as \emph{meta-symplectic}, quickly reviewed below. It is worth stressing that our main gadget, i.e., the correspondence between $\E$ and $\M^\E$, is just an adaptation of similar techniques traditionally found in other areas of modern mathematics\footnote{B\"acklund transformation, double f\/ibration transform, Penrose transform, etc. (see also~\cite{138544} on this concern).}.
In consequence, in this paper there coexist dif\/ferential and algebraic geometric constructions, and the authors are aware that this may give a certain feeling of incompatibility. Section~\ref{sezChiarificatrice} below should help in this concern.

\subsection*{Structure of the paper}
Section \ref{SecPrelim} reviews the classical notions needed for the formulation of the main Theo\-rems~\ref{LemmaNotturnoQuasiMattinieroPromossoATeorema},~\ref{thEgr} and~\ref{thFormMag}, which is given in Section \ref{subDescrMainRes}, and the Sections \ref{secCono}, \ref{SecChar3ordPDEs} introduce new ideas, namely the \emph{three-fold orthogonality} and its relationship with the \emph{characteristics variety}, necessary to the proof of the main results, worked out in Section~\ref{secCentrale}.
In Section~\ref{secIntInt} we study the intermediate integrals of Goursat-type $3\Rd$ order MAEs in terms of their characteristics.

\subsection*{Notation and conventions}
All objects and morphisms herewith considered are assumed of class $C^\infty$, up to mild exceptions (see Section \ref{sezChiarificatrice} below).

We use the symbol ``$\P$'' in order to avoid too many repetitions of the sentence ``up to a~conformal factor''.
One-dimensional linear object may be identif\/ied with their generators. As a~rule, if $P\to X$ is a bundle, and $x\in X$, we denote by $P_x$ the f\/ibre of it at~$x$. Symbol $\D^\prime$ denotes the derived distribution $\D+[\D,\D]$ of a distribution $\D$. Dif\/ferential forms on a mani\-fold~$N$ (resp., distribution $\D$) are denoted by~$\Lambda^\ast N$ (resp.,~$\Lambda^\ast \D^\ast$) and~$f_\ast$ denotes the tangent map of~$f$.
If $T$ is a tensor on $N$, or a distribution, we sometimes skip the index~``$n$'' in~$T_n$, if it is clear from the context that~$T$ has been evaluated in $n\in N$.
 The symmetric tensor product is denoted by~$\odot$.
The Einstein summation convention will be used, unless otherwise specif\/ied. Also, when a pair~$jk$ (resp., a~triple~$ijk$) runs in a~summation, such summation is performed over $1\leq j\leq k\leq 2$ (resp., $1\leq i\leq j\leq k\leq 2$), unless stated dif\/ferently.
Throughout the paper, by~$M^{(0)}$,~$\C^0$, and~$m^{(0)}$ we shall always mean~$M$, $\C$, and~$m$, respectively.

\begin{Remark}\label{remVerticale}
 If there is any bundle structure, or a vertical distribution, in the surrounding manifold, the symbol~$\D^v$ denotes the vertical part of the distribution~$\D$.
\end{Remark}

\section[Preliminaries on (prolongations of) contact manifolds, and (meta)symplectic structures]{Preliminaries on (prolongations of) contact manifolds,\\ and (meta)symplectic structures}\label{SecPrelim}

\subsection{Contact manifolds, their prolongations and PDEs}\label{sec.definitions.contact.and.prol}

Throughout this paper, $(M,\C)$ will be a 5D contact manifold, i.e., $\C$ is a completely non-integrable distribution of hyperplanes on $M$ locally described as $\C=\ker\theta$, where the $1$-form~$\theta$ is determined up to a~conformal factor and $\theta\wedge d\theta\wedge d\theta\neq 0$. The restriction $d \theta|_{\mathcal{C}}$
def\/ines on each~$\mathcal{C}_m$, $m\in M$, a conformal symplectic structure: Lagrangian (i.e., maximally $d\theta$-isotropic) planes of~$\mathcal{C}_m$ are tangent to maximal integral submanifolds of~$\mathcal{C}$ and, as such, their dimension is~$2$. We denote by $\mathcal{L}(\mathcal{C}_m)$ the \emph{Grassmannian of Lagrangian planes} of~$\mathcal{C}_m$ and by
\begin{gather}\label{eqDefEmmeUno}
\pi\colon \ M^{(1)} := \coprod_{m \in M}\mathcal{L}(\mathcal{C}_m) \to M
\end{gather}
the bundle of Lagrangian planes, also known as the $1\St$ \emph{prolongation} of $M$. The key property of the manifold $M^{(1)}$ is that it is naturally endowed
with a 5D distribution, def\/ined by
\begin{gather*}
\mathcal{C}^{1}_{m^1}:=\big\{\w\in T_{m^1}M^{(1)}\,|\,
\pi_{*}(\w)\in L_{m^1}\big\},\end{gather*}
where $L_{m^1}\equiv m^1$ is a point of $M^{(1)}$ considered as
a Lagrangian plane in $\mathcal{C}_{m}$.
Let us denote by~$\theta^{(1)}$ the set of $1$-forms on $M^{(1)}$ vanishing on $\mathcal{C}^{1}$. Then, by def\/inition, a \emph{Lagrangian plane} of $M^{(1)}$ is a~2D subspace which is $\pi$-horizontal\footnote{If horizontality is dropped in \eqref{eqDefEmmeDue}, one augments $M^{(2)}$ with the so-called ``singular'' integral elements of $\C^1$ which, in the case of PDEs, formalise the notion of ``singular solutions'', originally introduced in \cite{MR966202} (see also the review paper \cite{LucaChar} and references therein).} and such that all the forms belonging to the dif\/ferential ideal generated by $\theta^{(1)}$, vanish on it. In analogy with \eqref{eqDefEmmeUno}, we def\/ine the $2\Nd$ \emph{prolongation} $M^{(2)}$ of a contact
manifold $(M, \mathcal{C})$ as the f\/irst prolongation of $M^{(1)}$, that is
\begin{gather}\label{eqDefEmmeDue}
M^{(2)}=\big(M^{(1)}\big)^{(1)}:=\big\{ \text{Lagrangian planes of $M^{(1)}$} \big\}.
\end{gather}
Projection \eqref{eqDefEmmeUno} is the beginning of a tower of natural bundles $M^{(2)} \overset{\pi_{2,1}}{\longrightarrow} M^{(1)} \overset{\pi}{\longrightarrow} M
$,
which, for the present purposes, will be exploited only up to its $2\Nd$ term.
It is well known that $\pi_{2,1}$ is an af\/f\/ine bundle (see, e.g., \cite{MR1202431}).

The \emph{tautological bundle} $L\rightarrow M^{(i)}$
is def\/ined by requiring that the f\/ibre $L_{m^i}$ is $m^i$ itself, understood as a 2D subspace of $\C^{i-1}$, with $i=1,2$. We keep the same symbol $L$ for both the tautological bundles over $M^{(1)}$ and $M^{(2)}$, since it will be clear from the context which is which.

A generic point of $M$ (resp., $M ^{(1)}$, $M ^{(2)}$) is denoted by $m$ (resp., $m^1$, $m^2$). As a rule, when both $m^1$ and $m^2$ (resp. $m$ and $m^1$) appear in the same context, the former is always the $\pi_{2,1}$-image (resp., the $\pi$-image) of the latter.

By a $k\Th$ order PDE we always mean a sub-bundle $\E\subseteq M^{(k-1)}$ of codimension one whose f\/ibre $\E_{m^{k-2}}$ at $m^{k-2}$ is henceforth assumed, without loss of generality, to be connected (see~\eqref{sostituzionePerIGettofili} below). Recalling that $\pi_{2,1}$ is an af\/f\/ine bundle, we say that an equation $\E$ is \emph{quasi-linear} at a~point~$m^1\in M^{(1)}$ if the f\/ibre~$\E_{m^1}$ is an af\/f\/ine subspace of $\pi_{2,1}^{-1}(m^1)$, otherwise we say that it is \emph{non-linear} at $m^1$.
We retain the same symbol~$L$ for the tautological bundle $L|_\E\longrightarrow \E$ restricted to~$\E$.

\subsection[The meta-symplectic structure on $\C^1$]{The meta-symplectic structure on $\boldsymbol{\C^1}$}\label{subMetaSim}

The (local) conformal symplectic structure $\varpi=d\theta|_{\C}$ admits a ``global'' analog, namely
\begin{gather}\label{formaSimpletticaTwisted}
\varpi_{\textrm{glob}}\colon \ \C\wedge\C\longrightarrow \frac{\C^\prime}{\C}=\frac{TM}{\C}.
\end{gather}
Indeed,
since
$\frac{TM}{\C} $ is rank-one, $\varpi_{\textrm{glob}}$ locally identif\/ies with $\varpi$.

One of the main gadgets of our analysis is the co-restriction to $(\C^{1})^\prime$ of the Levi form of $\C^1$, f\/irstly investigated by V. Lychagin \cite{MR2389645} as a ``twisted'' analog of the symplectic form \eqref{formaSimpletticaTwisted} for the prolonged contact distribution $\C^1$ and called, for this reason, \emph{meta-symplectic}:
\begin{gather}\label{eqFormMetaSymp}
\Omega\colon \ \C^{1}\wedge\C^{1}\longrightarrow \frac{(\C^{1})^\prime}{\C^{1}}.
\end{gather}
Notice that,
unlike~\eqref{formaSimpletticaTwisted}, the form~\eqref{eqFormMetaSymp} takes its values into a rank-two bundle so that, even locally, it cannot be regarded as a 2-form in the standard sense. Nevertheless, it can be used for def\/ining Lagrangian subspaces: indeed, a~$\pi_{2,1}$-horizontal 2D subspace of $\C_{m^1}^1$ is \emph{Lagrangian} if it is $\Omega$-isotropic.

\section{Description of the main results}\label{subDescrMainRes}

In the case of a $3\Rd$ order PDE $\E$, polynomial \eqref{characteristic} always admits a linear factor, so that the corresponding characteristic cone $\M^\E$, \emph{always} contains a linear irreducible component $\MM^\E$ (we stress that in the fully decomposable case, all the three irreducible components are linear, and it suf\/f\/ices to choose one).

\begin{Definition}
 If $\omega\in \Lambda^2 M^{(1)}$ is a 2-form, then the hypersurface
 \begin{gather}\label{eqEOmega}
\E_{\omega}:=\big\{ m^2\in M^{(2)}\,|\, \omega|_{L_{m^2}}= 0 \big\}
\end{gather}
in $M^{(2)}$ is called the \emph{Boillat-type} $3\Rd$ order MAE (associated to $\omega$).
\end{Definition}

The f\/irst result of this paper is concerned with the characteristic cone $\M^{\E_\omega}$. 
 Namely, we show that $\M^{\E_\omega}$ can be used to recover the equation itself, since there is an obvious way to ``invert'' the construction of the characteristic cone \eqref{eq.char.cone.1} out of a PDE $\E$: take \emph{any} sub-bundle $\M\subseteq\P\C^1$ and associate with it the following subset
\begin{gather}\label{eqENu}
\E_\M:=\big\{ m^2 \in M^{(2)}\,|\, \exists\,\elle\in\M\colon\, L_{m^2}\supset\elle\}\subseteq M^{(2)}.
\end{gather}
Now, if $\M$ is regular enough, the corresponding $\E_\M$ turns out to be a genuine $3\Rd$ order PDE, and examples of ``regular enough'' sub-bundles are provided by characteristic cones of~$3\Rd$ order PDEs themselves. In particular, it always holds the inclusion $\E_{\M^\E}\supseteq \E$.

\begin{Definition}\label{defMaledetta}
 We say that the equation $\E$ is \emph{recoverable} (from its characteristics) if
 \begin{gather}
\E=\E_{\M^\E}.\label{eqStar}
\end{gather}
\end{Definition}

\begin{Theorem} \label{LemmaNotturnoQuasiMattinieroPromossoATeorema}
 Any $3\Rd$ order MAE is recoverable from its characteristics.
 \end{Theorem}

The simplest examples of equations \eqref{eqENu} are obtained when $\M=\P\D$, where $\D \subset \C^1$ is a 3D sub-distribution:
 \begin{gather}\label{eqED}
\E_{\D }:=\E_{\P\D }=\big\{ m^2\in M^{(2)}\,|\, L_{m^2}\cap\D_{m^1} \neq 0 \big\},
\end{gather}
with $m^2$ projecting onto $m^1$. Equations~\eqref{eqED} are herewith dubbed \emph{Goursat-type~$3\Rd$ order MAEs}.
Observe that the equations \eqref{eqED} form the sub-class of the equations~\eqref{eqEOmega} which are determined by 2-forms which are decomposable modulo the dif\/ferential ideal generated by contact forms (see the beginning of Section~\ref{subProofThEgr}).
The second result of this paper allows to locally characterizes Goursat-type $3\Rd$ order MAEs in terms of their characteristic cone.

\begin{Theorem}\label{thEgr}
Let $\E$ be a $3\Rd$ order PDE. Then $\E$ is locally of the form \eqref{eqED}, if and only if its characteristic cone $\M^\E$ contains an irreducible component $\MM^\E$ which is a $2D$ linear projective sub-bundle.
\end{Theorem}

Theorem \ref{thEgr}
implies that $\E\subseteq M^{(2)}$ is a Goursat-type $3\Rd$ order MAE if and only if its characteristic cone decomposes as
\begin{gather}\label{decCono}
 \M^\E= \P\D_1\cup\MMM^\E,
\end{gather}
where $\D_1\subseteq \C^1$ is a 3D sub-distribution and $\MMM^\E$ encompasses all the remaining irreducible components of~$\M^\E$. Surprisingly enough, if \eqref{decCono} holds, then all the information about~$\E$ is encapsulated in~$\D_1$.
For instance, if $\MMM^\E$ is in its turn reducible, then its components~$\P\D_2$ and~$\P\D_3$ are linear as well and can be unambiguously characterized by $\P\D_1$ through the formula
 \begin{gather}\label{eqFormulaMagica}
\Omega\big(h_1,\D_i^v\big)=\Omega\big(h_i,\D_1^v\big)=\textrm{1D space},
\end{gather}
where $h_i$ is a non-zero horizontal vector in $\D_i$, $i=1,2,3$ (for the def\/inition of $\D^v$ see Remark~\ref{remVerticale}).

Observe that, if a Goursat-type MAE is non-linear in a point $m^1\in M^{(1)}$, then it is non-linear in a neighbourhood of~$m^1$; on the contrary, if it is linear in $m^1$, then a neighbourhood of $m^1$ where it is quasi-linear may not exist (see~\eqref{eqEsempioType}). On account of this, we give the following Theorem~\ref{thFormMag}, that f\/inalises our characterisation of Goursat-type $3\Rd$ order MAEs through their characteristics.

\begin{Theorem}
\label{thFormMag}
 Let $\E=\E_{\P\D_1}$, where $\D_1\subset\C$ is a $3$-dimensional distribution, be a Goursat-type MAE. Let $m^1\in M^{(1)}$. Then $\E=\E_{\P\D_1\cup\MMM^\E}$ and
\begin{enumerate}\itemsep=0pt
 \item[$1)$] $\E$ is quasi-linear at $m^1$ if and only if $\dim(\D_1^v)_{m^1}=2$: in this case, $(\MMM^\E)_{m^1}$ is either empty or equal to $\P({\D_2})_{m^1}\cup\P({\D_3})_{m^1}$, where $({\D_i})_{m^1}$ is unambiguously defined by \eqref{eqFormulaMagica} and $\E_{(\D_i)_{m^1}}=\E_{(\D_1)_{m^1}}$, $i=2,3$;
\item[$2)$] $\E$ is non-linear at $m^1$ if and only if $\dim(\D_1^v)_{m^1}=1$: in this case, if $(\MMM^\E)_{m^1}$ is not empty, then it cannot contain any linear irreducible component and $\E_{m^1}=\E_{(\MMM^\E)_{m^1}}$;
\item[$3)$] if $\E_{m^1}=\E_{(\D_2)_{m^1}}$ for some $3D$ subspace $(\D_2)_{m^1}\subset\C^1_{m^1}$, then either $(\D_1)_{m^1}=(\D_2)_{m^1}$ or $(\D_1)_{m^1}$ and $(\D_2)_{m^1}$ are ``orthogonal'' in the sense of~\eqref{eqFormulaMagica}.
\end{enumerate}
\end{Theorem}

As a main byproduct of Theorem~\ref{thFormMag} we shall obtain a method for f\/inding intermediate integrals, discussed in Section~\ref{secIntInt}.
We conclude this section with a local description of the main objects introduced so far.

\subsection{Local coordinate description of the main objects}

Let $(x^i, u, p_i)$ be contact coordinates on $M$, $i=1, 2$. Then
$\theta = du - p_i dx^i$, and
\begin{gather}\label{eqCoordCi}
\C = \Span{D_1,D_2, \partial_{p_1},\partial_{p_2}},
\end{gather}
where $D_i$ is the total derivative with respect to $x^i$, truncated to the $0\Th$ order. The above-introduced system $(x^i,u,p_i)$ induces coordinates
\begin{gather}\label{coordinates_prolongation}
\big(x^i, u, p_i, p_{ij}= p_{ji} ,\, 1 \leq i, j \leq 2 \big)
\end{gather}
on $M^{(1)}$ as follows: a point $m^1 \equiv L_{m^1} \in M^{(1)}$
has coordinates \eqref{coordinates_prolongation} if\/f $m=\pi(m^1)=(x^i,u,p_i)$ and the
corresponding Lagrangian plane $L_{m^1}$ is given by
\begin{gather}\label{eqCoordEllm1}
L_{m^1}=\langle D_i + p_{ij}\partial_{p_{j}}\,| \,i\in\{1,2\}\rangle\subset\mathcal{C}_{m}.
\end{gather}
Similarly,
\begin{gather}\label{eqCoordCiUno}
\CC=\Span{D_1,D_2, \partial_{p_{11}},\partial_{p_{12}},\partial_{p_{22}}},
\end{gather}
where now the $D_i$'s stand for the total derivatives truncated to the $1\St$ order and a point $m^2\in M^{(2)}$ has coordinates $(x^i,u,p_i,p_{ij}=p_{ji},p_{ijk}=p_{ikj}=p_{jik}=p_{jki}=p_{kij}=p_{kji})$ if the corresponding Lagrangian plane is given by
\begin{gather}\label{eqCoordEll}
L_{m^2}=\Span{ D_i+ p_{ijk}\partial_{p_{jk}}\,| \,i\in\{1,2\}}\subseteq \C^1_{m^1},
\end{gather}
where the vector f\/ields $D_i$ and $\partial_{p_{ij}}$ are tacitly assumed to be evaluated at $m^1$.

\begin{Remark}\label{remcoordELLE}
We always use the symbol $D_i$ for the truncated total derivative with respect to $x^i$, $i=1,2$, the order of truncation depending on the context. For instance, the order of truncation is~$0$ in~\eqref{eqCoordCi} and it is $1$ in~\eqref{eqCoordCiUno} and~\eqref{eqCoordEll}.
It is convenient to set
 \begin{gather}\label{defYi}
\xi_i:=\left. D_i\right|_{m^k},\qquad k=1,2,
\end{gather}
where the total derivatives appearing in~\eqref{defYi} are truncated to the $(k-1)\St$ order. Indeed, both~\eqref{eqCoordEllm1}
and~\eqref{eqCoordEll} simplify as
\begin{gather}
L_{m^k}=\Span{ \xi_1,\xi_2 },\qquad k=1,2,\label{coordELLE}
\end{gather}
and their dual as
$L_{m^k}^\ast=\Span{dx^1,dx^2}$,
respectively.
\end{Remark}

Using the local coordinates
\eqref{eqCoordCiUno}, it is easy to realize that $(\C^{1})^\prime$ is spanned by $\C^1$ and $\Span{\partial_{p_1},\partial_{p_2}}$, so that the quotient $\frac{(\C^{1})^\prime}{\C^1}$ identif\/ies with the latter, and \eqref{eqFormMetaSymp} reads
 \begin{gather}\label{eqFormMetaSympLOC}
\Omega=dp_{ij} \wedge dx^i \otimes\partial_{p_j}.
\end{gather}
As a vector-valued dif\/ferential 2-form, $\Omega$ can be identif\/ied with the pair $\Omega \equiv (-d\theta_1,-d\theta_2)$
where $\theta_i=dp_i-p_{ij}dx^j$.

\subsubsection[Boillat and Goursat $3\Rd$ order MAEs]{Boillat and Goursat $\boldsymbol{3\Rd}$ order MAEs}

We conclude with a local coordinate description of Boillat \eqref{eqEOmega} and Goursat \eqref{eqED} equations $\E=\{F=0\}$.
Locally, the former is given by
\begin{gather}\label{eq.general.MAE.3}
F=\det\left(\begin{matrix}p_{111} & p_{112} & p_{122} \\p_{112} & p_{122} & p_{222} \\ & \A & \end{matrix}\right)+\B\cdot (p_{111},p_{112},p_{122},p_{222})^T+C,
\end{gather}
where $ \A$ (resp., $\B$, $C$) is an $\R^3$-(resp., $\R^4$-, $\R$-)valued smooth function on $M^{(1)}$,
and the latter is either given by
 \begin{gather}\label{eqGoursat}
 F=\det\left(\begin{matrix}p_{111}-f_{111} & p_{112}-f_{112} & p_{122}-f_{122} \\p_{112}-f_{211} & p_{122}-f_{212} & p_{222}-f_{222} \\ & \A & \end{matrix}\right),
\end{gather}
where $\A$ is as in \eqref{eq.general.MAE.3} and $f_{ijk}\in C^\infty(M^{(1)})$, or it is quasi-linear.

\section[Vertical geometry of (prolongations of) contact manifolds and their characteristics]{Vertical geometry of (prolongations of) contact manifolds\\ and their characteristics}\label{secCono}

The departing point of our analysis of $3\Rd$ order PDEs is to identify them with sub-bundles of the~$1\St$ prolongation $M^{(2)}=(M^{(1)})^{(1)}$ of~$M^{(1)}$. Hence, a key role will be played by their vertical geometry, i.e., the $1\St$ order approximation of their bundle structure, and, in particular, by the so-called \emph{rank-one vectors}, which are in turn linked to the notion of \emph{characteristics}. In order to introduce these concepts, we begin with the vertical geometry of the surrounding bundle, i.e.,~$M^{(2)}$ itself.

\subsection[Vertical geometry of $M^{(k)}$ and three-fold orthogonality in $M^{(1)}$]{Vertical geometry of $\boldsymbol{M^{(k)}}$ and three-fold orthogonality in $\boldsymbol{M^{(1)}}$}

 For $k\in\{1,2\}$,
we def\/ine the vertical bundle over $M^{(k)}$ as follows
\begin{gather*}
VM^{(k)}:=\coprod_{m^k\in M^{(k)}}T_{m^k}M^{(k)}_{m^{k-1}}.
\end{gather*}
As regard to the contact manifold $M$, i.e., the case $k=0$, it does not possess a naturally def\/ined vertical bundle: we replace it by the following bundle on $M^{(1)}$:
\begin{gather}\label{eqDefRelVecBund}
VM:=\coprod_{m^1\in M^{(1)}}\frac{\C_{\pi (m^1)}}{L_{m^1}}.
\end{gather}

\begin{Lemma}\label{lemFund}
It holds the following canonical isomorphism:
 \begin{gather}\label{isoFond}
VM^{(k)}\simeq S^{k+1}L^*, \qquad k\geq 0.
\end{gather}
\end{Lemma}
\begin{proof}
It is a generalisation of the proof of the jet-theoretic version of the statement (see, e.g., \cite[Theo\-rem~3.2]{ThesisMichi}).
\end{proof}

Directly from the def\/inition \eqref{eqDefRelVecBund} of $VM$ and the isomorphism \eqref{isoFond} for $k=0$ one obtains
\begin{gather}\label{eq1ellstar}
 \frac{\C_{\pi (m^1)}}{L_{m^1}}\cong L_{m^1}^\ast.
\end{gather}

\begin{Example}
Of particular importance will be the vertical vectors on $M^{(2)}$ which correspond to the \emph{perfect cubes} of covectors on $L$ via the fundamental isomorphism \eqref{isoFond}. For instance, if
\begin{gather}\label{EqCoordAlpha}
\alpha:=\alpha_idx^i\in L_{m^2}^\ast
\end{gather}
(see Remark \ref{remcoordELLE}), then the corresponding vertical vector is
\begin{gather}\label{isoFondCoord}
S^3L_{m^2}^\ast\ni\alpha^{\otimes 3}\longleftrightarrow \sum_{i+j=3} \alpha_1^i\alpha_2^j \frac{\partial }{\partial p_{\underbrace{\text{\scriptsize $1\cdots 1$}}_{i} \underbrace{\text{\scriptsize $2\cdots 2$}}_{j}}} \Bigg|_{ {m}^{2}}\in V_{m^2}M^{(2)}.
\end{gather}
\end{Example}

Even if the sections of $VM$ are \emph{not}, strictly speaking, vector f\/ields, and a such they lack an immediate geometric interpretation, the bundle~$VM$ itself has important relationships with the contact bundle~$\C^1$. Namely, on one hand, it is canonically embedded into the module of 1-forms on~$\C^1$ (see~\eqref{StranissimoEmbedding} later on) and, on the other hand, it is identif\/ied with the quotient distribution~$\frac{(\C^1)^\prime}{\C^1}$ (Lemma~\ref{lemmaStupidoMaCarino} below).

\begin{Lemma}\label{lemmaStupidoMaCarino}
There is a $($conformal$)$\footnote{In the sense that \eqref{eq21ellstar} is canonical only up to a conformal factor.} natural isomorphism
\begin{gather}\label{eq21ellstar}
\frac{(\C^1)^\prime}{\C^1}\cong L^\ast.
\end{gather}
\end{Lemma}

\begin{proof}
It follows from a natural isomorphsim between the left-hand sides of \eqref{eq1ellstar} and \eqref{eq21ellstar}. Indeed, the map
\begin{gather}\label{eqInstrinsic}
\frac{(\C^1)^\prime_{m^1}}{\C^1_{m^1}}\longrightarrow \frac{\C_{\pi(m^1)}}{L_{m^1}}
\end{gather}
induced from $\pi_\ast$ is well-def\/ined and linear, for all $m^1\in M^{(1)}$. In other words, \eqref{eqInstrinsic} is a bundle morphism, is well-def\/ined and linear, and surjective. Since the ranks of the bundles are the same, it is an isomorphism.
\end{proof}

We stress that there is no canonical way to project the bundle $\C^1$ over the tautological bundle~$L$, if both are understood as bundles over $M^{(2)}$. Nevertheless, if $\C^1$ and $L$ are regarded as bundles over~$M^{(1)}$, then $\pi_*$ turns out to be a bundle epimorphism from $\C^1$ to $L$, i.e.,
 \begin{gather}\label{eqOrizzontalizzazione}
 \pi_*\colon\ \C^1\longrightarrow L=\pi_*\big(\C^1\big).
\end{gather}

Dually, epimorphism \eqref{eqOrizzontalizzazione} leads to the bundle embedding $L^\ast \hookrightarrow \C^{1\,\ast}$ which can be combined with the identif\/ication $L^\ast\cong VM$. The result is a (conformal) embedding of bundles over~$M^{(1)}$,
\begin{gather}\label{StranissimoEmbedding}
VM\hookrightarrow \C^{1\,\ast},
\end{gather}
which will be useful in the sequel.
In local coordinates, \eqref{StranissimoEmbedding} reads
$
\left. \partial_{p_i} \right|_{m}\longmapsto d_{m^1}x^i$, $ i=1,2
$.

Now we are in position to def\/ine the concept of orthogonality in the meta-symplectic context
(that has apparently never been observed before), which
generalizes the symplectic orthogonality within the contact distribution $\C$ of $M$. Indeed, an immediate consequence of Lemma \ref{lemmaStupidoMaCarino} is that
the meta-symplectic form $\Omega$ is $L^\ast$-valued, i.e., $\C^1\wedge\C^1\stackrel{\Omega}{\longrightarrow}L^\ast$ is a trilinear form
\begin{gather}\label{eqFormaProtoTrilineare}
\C^1\wedge\C^1\otimes_{M^{(1)}}L\stackrel{\Omega}{\longrightarrow}C^\infty\big(M^{(1)}\big).
\end{gather}
In turn, thanks to
the canonical projection \eqref{eqOrizzontalizzazione} of $\C^1$ over $L$, the form \eqref{eqFormaProtoTrilineare} descends to a
trilinear form $\C^1\wedge\C^1\otimes \C^1\stackrel{\widetilde{\Omega}}{\longrightarrow}C^\infty(M^{(1)})$.
Such notions (orthogonal vectors, orthogonal complement, Lagrangian subspaces) can be found also in a 5D meta-symplectic space, but with more subtleties. For example, there can be up to \emph{two} distinct, so to speak, ``orthogonal complements'' to a given subspace.

\begin{Definition}[three-fold orthogonality]\label{defTriOrt}
Elements $X_1,X_2,X_3\in \C^1$ are \emph{orthogonal} if and only if
$\widetilde{\Omega}(X_1,X_2,X_3)=0$.
\end{Definition}

\begin{Remark}\label{remFomulaMagicaIntrinseca}
In this perspective, condition \eqref{eqFormulaMagica} expresses precisely the fact that the pair $\{\D_2,\D_3\}$ is the orthogonal complement of the distribution $\D_1$, in the sense that
\begin{gather*}
\widetilde{\Omega}(\D_i,\D_j,\D_k)=0,\qquad \{i,j,k\}=\{1,2,3\}.
\end{gather*}
\end{Remark}

\subsection{Rank-one lines, characteristic directions and characteristic hyperplanes}\label{subsubRankUnoChar}

Isomorphism \eqref{isoFond} shows that there is a canonical distinguished subset of tangent directions to $M^{(k)}_{m^{k-1}}$, $k\in\{1,2\}$, namely those sitting in the image of the Veronese embedding $\P L^\ast\hookrightarrow \P S^{k+1}L^\ast$, i.e., the $k\Th$ powers of sections of $L^\ast $. Geometrically, these are the tangent directions at $m^k=m^k(0)$ to the curves $m^k(t)$ such that the corresponding family $L_{m^k(t)}$ of Lagrangian subspaces in $M^{(k-1)}$ ``rotates'' around a common hyperplane (which, in our case, is a line).

Local sections of the dual bundle $\C^{k\ast}$ will be called \emph{vertical forms}, as they are the image of the local sections of $V^*M^{(k)}$ under the natural projection\footnote{We insist that ``vertical forms'' are \emph{not}, globally speaking, dif\/ferential forms on $M^{(k)}$, but rather equivalence classes of them. Locally, of course, one can f\/ind a representative in each class.}. Because of the dual of the fundamental isomorphism (see Lemma~\ref{lemFund}), $\C^{k\ast}$ is the epimorphic image of the symmetric power~$S^{k+1}L$ of the tautological bundle, and we can speak of \emph{decomposable vertical forms}, if they come from exact powers.

\begin{Definition}\label{defRankOneLine}
 A line $\ell\in \mathbb{P}V_{m^k}M^{(k)}$ is called a \emph{rank-one line} if $\ell=\Span{\alpha^{\otimes (k+1)}}$, for some $\Span{\alpha}\in \P L_{m^k}^\ast$, in which case we call $\alpha$ a \emph{characteristic covector} and $\Span{\alpha}$ a \emph{characteristic $($direction$)$} in the point $m^k$. The subspace $H_\alpha:=\ker\alpha\leq L_{m^k}$ is called the \emph{characteristic hyperplane} associated to the rank-one line $\ell$. Furthermore, if $\omega\in\C^{k\,\ast}$ is a vertical form, a hyperplane $H_\alpha$ is called \emph{characteristic} for $\omega$ if $\omega|_\ell=0$.
\end{Definition}

Denomination \emph{rank-one} refers to the rank of the multi-linear symmetric form on $L$ involved in the def\/inition (see, e.g., \cite{MR2369200}).
Observe that, in our contest, $\dim H_\alpha=1$, so that we shall speak of a \emph{characteristic line} and call \emph{characteristic vector} a generator of $H_{\alpha}$. The geometric relationship between $H_\alpha$ and $\ell=\Span{\alpha^{\otimes 3}}$ is well-known (see, e.g., \cite{LucaChar}) and it can be rendered by
\begin{gather}\label{eqGeometriaDelRaggio}
\ell=\P T_{m^k}\big( H_\alpha^{(1)} \big),
\end{gather}
where the $1\St$ {prolongation} $ H_\alpha^{(1)}=\{ {m}^k\in M^{(k)}_{m^{k-1}}\,|\, L_{ {m}^k}\supseteq H_\alpha \}$ of $H_\alpha$ has dimension one.

By \emph{prolongation} of a subspace $W\subseteq T_{m^k}M^{(k)}$ we mean the subset $W^{(1)}\subseteq M^{(k+1)}_{m^k}$ made of the Lagrangian planes which are contained in $W$. The prolongation of a submanifold $\mathcal{W}\subseteq M^{(k)}$ is the sub-bundle of $\mathcal{W}^{(1)}\subseteq M^{(k+1)}$ def\/ined by $\mathcal{W}^{(1)}_{m^k}:=(T_{m^k}\mathcal{W})^{(1)}$, $m^k\in M^{(k)}$. The prolongation to $M^{(k)}$ of a local contactomorphism $\psi$ of $M$ is denoted by $\psi^{(k)}$.

\begin{Remark}\label{RemarcoFondamentale}
 Let $\alpha$ and $L_{m^2}$ be locally described by \eqref{EqCoordAlpha} and \eqref{defYi}. Then $ H_\alpha= \ker(\alpha_1dx^1+\alpha_2dx^2)=\Span{\alpha_2 \xi_1,-\alpha_1 \xi _2}$,
i.e., one can identify $\alpha\leftrightarrow H_\alpha$ via a counterclockwise $\frac{\pi}{2}$ rotation:
\begin{gather}\label{RotazioneDelCazzo}
\alpha\equiv (\alpha_1,\alpha_2)\leftrightarrow(\alpha_2,-\alpha_1)\equiv H_\alpha.
\end{gather}
From now on, characteristic directions and characteristic lines
 will be taken as synonyms, thus breaking the separation proposed below:
\begin{center}
\begin{tabular}{c|c|c}
&contravariant object& covariant object \\
\hline
$M^{(2)}$&$\underset{\textrm{rank-one line}}{\ell=T_{m^2} H_\alpha^{(1)}}$ & \\
\hline
$M^{(1)}$&$\underset{\textrm{characteristic hyperplane}\equiv\textrm{line}}{H_\alpha=\ker\alpha\in \P L_{m^2} }$ & $\underset{\textrm{characteristic direction}}{\Span{\alpha}\in\P L_{m^2}^\ast,}\ \underset{\textrm{characteristic covector}}{\alpha\in L_{m^2}^\ast{\setminus}\{0\}}$
\end{tabular}
\end{center}

It should be stressed that in the multidimensional cases (see, e.g.,~\cite{MR2985508})
it is no longer possible to regard the characteristic directions as lines lying in~$L$, since the latter are \emph{not} hyperplanes.
\end{Remark}

\subsection{Canonical directions associated with orthogonal distributions}

Suppose that the vertical form $\omega\in \C^{1\,\ast}$ is decomposable in the sense of Lemma \ref{lemFund}. Then it is easy to see that there are two characteristic lines $\elle_1$ and $\elle_2$.
The meta-symplectic form \eqref{eqEOmega} links these two lines with the vertical distribution $V:=\ker\omega$ determined by $\omega$.

Recall now that $\Omega$ is $VM\simeq L^*$-valued (see \eqref{isoFond} and \eqref{eqFormaProtoTrilineare}), so that $H_i$ can be considered as a subspace of $VM$ (in view of Remark \ref{RemarcoFondamentale}), and there are the induced maps
\begin{gather*}
 \Omega(X,\,\cdot\, )\colon \ \C^1\longrightarrow \frac{VM}{H_i} ,\qquad X\in\C^1 .
\end{gather*}
Letting $X$ vary in $V$, one gets the rank$\leq 2$ sub-bundle
\begin{gather*}
\frac{V\lrcorner \Omega}{H_i}:=\Span{ \Omega(X,\,\cdot\, )+ H_i\,|\, X\in V }\subset \C^{1\ast}\otimes \frac{VM}{H_i} .
\end{gather*}
Finally, $\elle_j\otimes\frac{L^\ast}{\elle_i}$ can be regarded as a one-dimensional sub-bundle of $\C^{1\ast}\otimes \frac{VM}{H_i}$, via \eqref{StranissimoEmbedding}.

Taking this into account, we have the following result.

\begin{Lemma}\label{lemma.strano}It holds the following identification
\begin{gather*}
 \frac{V\lrcorner \Omega}{\elle_i}\simeq\elle_j\otimes\frac{L^\ast}{\elle_i},\qquad \{i,j\}=\{1,2\}.
\end{gather*}
In particular,
$\dim\frac{V\lrcorner \Omega}{\elle_i}=1$, $ i=1,2$.
\end{Lemma}

\begin{proof}
We shall assume that the coef\/f\/icient of $dp_{11}$ is not zero: other cases can be dealt with likewise. So, being $\omega$ decomposable, it can be brought, up to a scaling, to the form
$\omega=dp_{11}-(k_1+k_2)dp_{12}+k_1k_2dp_{22}\in\C^{1\ast}$.
Accordingly,
$\elle_i=\partial_{p_2}+k_i\partial_{p_1}$
and a straightforward computation shows that
$V=\Span{X_1,X_2}$,
with
\begin{gather}
X_1=-k_1k_2\partial_{p_{11}}+\partial_{p_{22}},\qquad X_2=(k_1+k_2)\partial_{p_{11}}+\partial_{p_{12}}.\label{eqX1X2}
\end{gather}
Then
\begin{align}
V\lrcorner \Omega &=\Span{X_1\lrcorner \Omega,X_2\lrcorner \Omega} \nonumber \\
 &
 =\Span{k_1k_2 dx^1\otimes \partial_{p_1}-dx^2\otimes \partial_{p_2}, -(k_1+k_2)dx^1\otimes \partial_{p_1}-dx^1\otimes \partial_{p_2}-dx^2\otimes \partial_{p_1}} .\label{VInserOmega}
\end{align}
If the basis $\{\elle_i,\partial_{p_1}\}$ for $VM\simeq L^\ast$ (see \eqref{isoFond}) is chosen, then the factor of \eqref{VInserOmega} with respect to $\elle_i$ is the line
\begin{gather}\label{eqLineaFinale}
\Span{\big(k_jdx^1+dx^2\big)\otimes\partial_{p_1}},\qquad i\neq j,
\end{gather}
and the lemma is proved.
\end{proof}

Lemma \ref{lemma.strano} will be used later on in the proof of Proposition~\ref{propAnnullatoreEDirezioni}, establishing that,
associated with a~quasi-linear MAE with completely decomposable symbol, there are three canonical lines in~$M$.

\subsection{Some examples}\label{exEqQL}

Let $\E$ be as in \eqref{equation}, with $F$ given by $p_{111}-p_{112}-2p_{122}$ (resp., $p_{122}$ and~$p_{111}$). Then equation~\eqref{characteristic} reads
\begin{gather}\label{EXcharacteristic}
\big(\w^2\big)^3+\big(\w^2\big)^2\w^1-2\w^2\big(\w^1\big)^2= 0\qquad \big(\textrm{resp.,} \ \w^2\big(\w^1\big)^2=0 \ \textrm{and}\ \big(\w^2\big)^3=0\big).
\end{gather}
Applying def\/inition~\eqref{eq.char.cone.1}, one obtains that $\M^\E=\cup_{i=1}^3\P\D_i$, with
\begin{gather}
\mathcal{D}_1 = \Span{ D_1 , \partial_{p_{11}}+\partial_{p_{12}}, 2\partial_{p_{11}}+\partial_{p_{22}}},\label{eq.piano.1}\\
\mathcal{D}_2 = \Span{ D_1 + D_2 ,2\partial_{p_{11}}+\partial_{p_{12}}, \partial_{p_{22}}},\label{eq.piano.2}\\
\mathcal{D}_3 = \Span{ D_1 - 2D_2 ,\partial_{p_{11}}-\partial_{p_{12}}, \partial_{p_{22}}}\label{eq.piano.3}
\end{gather}
(resp.,
 \begin{gather}
\mathcal{D}_1=\Span{ D_1,\partial_{p_{11}},\partial_{p_{12}}}, \qquad \mathcal{D}_2=\mathcal{D}_3=\Span{ D_2,\partial_{p_{11}},\partial_{p_{22}} }\label{eqServeAllaFineAlloraNonTogliere}
\end{gather}
and \looseness=-1
$
\mathcal{D}_1=\mathcal{D}_2=\mathcal{D}_3=\Span{ D_1,\,\partial_{p_{12}},\,\partial_{p_{22}} }
$).
Observe that, in the f\/irst case, equation \eqref{EXcharacteristic} factors as
\begin{gather}\label{EXcharacteristicSPLIT}
\w^2\big(\w^2-\w^1\big)\big(\w^2+2\w^1\big)=0,
\end{gather}
and that the ``horizontal components'' $h_i$ (i.e., the f\/irst generators) of~\eqref{eq.piano.1},~\eqref{eq.piano.2} and~\eqref{eq.piano.3} are precisely the three distinct roots of~\eqref{EXcharacteristicSPLIT}. Notice also that,
 in general, the distributions associated with the same quasi-linear $3\Rd$ order MAE need not to be contactomorphic. Here it follows an easy counterexample. For instance,
consider the equation $p_{122}=0$ above: its characteristic cone consists of the two distributions $\D_1$ and $\D_2=\D_3$ (see~\eqref{eqServeAllaFineAlloraNonTogliere}),
which are not contactomorphic. Indeed, the derived distribution $\D_1':=\D_1+[\D_1,\D_1]$ of $\D_1$ is $5$-dimensional, whereas $\dim(\D_2')=4$.

Let $\E$ be the f\/irst equation from above, i.e., with $F=p_{111}-p_{112}-2p_{122}$. The following three lines
 \begin{align}
1\St \ \textrm{line:} \quad & \Omega\big(\D_1,\D_2)=\Span{ 2\partial_{p_1}+\partial_{p_2}},\\
2\Nd \ \textrm{line:} \quad & \Omega\big(\D_1,\D_3)=\Span{ \partial_{p_1}-\partial_{p_2}},\label{MagDir2}\\
3\Rd \ \textrm{line:} \quad & \Omega\big(\D_2,\D_3)=\Span{ \partial_{p_2}},\label{MagDir3}
\end{align}
 are canonically associated with the triple $(\D_1,\D_2,\D_3)$, i.e., with the equation $\E$. It is worth observing that the vertical part $\D_1^v$ is an integrable vertical distribution on $M^{(1)}$, so that, in this case, lines \eqref{MagDir2} and \eqref{MagDir3} are the characteristic lines of the family of equations $2p_{22}-p_{11}+p_{12}=k(x^1,x^2,u)$.
Nevertheless, integrability $\D_1^v$ is not indispensable for associating directions \eqref{MagDir2} and \eqref{MagDir3} with the distribution $\D_1$, since they occur as the characteristic lines
of the (not necessarily closed) vertical covector $2dp_{22}-dp_{11}+dp_{12}$,
which annihilates $\D_1^v$.

Now we give another example that clarif\/ies the type-changing phenomenon:
\begin{gather}\label{eqEsempioType}
p_{11}(p_{112}p_{122}-p_{111}p_{222})-p_{122}=0 .
\end{gather}
In fact, equation \eqref{eqEsempioType} is non-linear at the points with $p_{11}\neq 0$, whereas it is linear at the points with $p_{11}=0$. Equation \eqref{eqEsempioType} is of Goursat type and the corresponding $3$-dimensional distribution
\begin{gather*}
\mathcal{D}=\langle D_1,\,p_{11}D_2+\partial_{p_{11}},\,\partial_{p_{12}} \rangle
\end{gather*}
behaves accordingly with Theorem \ref{thFormMag}, i.e., the dimension of $\mathcal{D}^v$ is $1$ if $p_{11}\neq 0$ and it is $2$ if $p_{11}=0$. Moreover, at the points with $p_{11}=0$ we have three planes, whose two of them coincide, def\/ining the f\/ibre of equation \eqref{eqEsempioType} (see \eqref{eqServeAllaFineAlloraNonTogliere}).

\subsection{Smoothness and singularities issues}\label{sezChiarificatrice}

For generic PDEs, we just assume that $\E\longrightarrow M^{(1)}$ is the zero locus of a smooth function, without requiring the non-vanishing of its dif\/ferential. In all the statements and reasonings, we tacitly restrict ourselves to the open and dense subset $\E_{\textrm{reg}}$ of smooth points.

According to our def\/inition of a PDE $\E$ as a sub-bundle of $M^{(1)}$, the natural projection $\E\longrightarrow M^{(1)}$ is surjective,
 and that each f\/ibre $\E_{m^1}$ is a smooth submanifold of $M^{(2)}_{m^1}$, possibly with singularities.

For example, the singularities which can occur for $3\Rd$ Monge--Amp\`ere equations are of ``algebraic type'', i.e., $\E\longrightarrow M^{(1)}$ is a bundle of projective varieties. Indeed, each f\/ibre $M^{(2)}_{m^1}$ is naturally understood as the Grassmannian variety of isotropic elements with respect to the meta-symplectic structure which, in analogy with the standard Lagrangian Grassmannian, embeds into a suitable projective space via the Pl\"ucker embedding\footnote{The general theory of this would lead us beyond the scope of the present paper.}. More precisely, a straightforward computation shows that the subset of the singular points of $\E_\D$ is made by the $m^2$'s such that $L_{m^2}\subset\D_{m^1}$. Such a possibility occurs only when $\dim(L_{m^2}\cap\D_{m^1})=2$ for $m^2\in M^{(2)}$ (in the case when $L_{m^2}$ and $\D_{m^1}$ are not transversal, generically the dimension is~1, otherwise is generically~0).

Methodologically, this paper is at the crossroad between algebraic and dif\/ferential geometry. The geometric approach to non-linear PDEs is traditionally dif\/ferential, while recent developments revealed that many features of PDEs are genuinely algebraic-geometric. Even though a~rigorous algebraic-geometric approach to the present topics is feasible, the authors have opted for the traditional and~-- in a sense, easier~-- methods, based on elementary dif\/ferential geometry.

Summing up, all the objects are herewith assumed to be smooth, and all maps to be of class~$C^\infty$, except for:
\begin{itemize}\itemsep=0pt
\item the f\/ibres of $\E$, which are the zero loci of smooth functions (and, as such, may display singularities where the dif\/ferential vanishes);
\item the sub-bundles $\M$ of the projectivised contact bundle $\P\C^1$, which are smooth families of projective sub-varieties in each f\/ibre, i.e., locally isomorphic to the product of an open subset of $\M$ by a projective subvariety in $\P^4$.
\end{itemize}

If $\MM$ and $\MMM$ are sub-bundles of $\M$, in the aforementioned sense, we say that $\MM$ and $\MMM$ are irreducible components of $\M$ if, for any point $m^1\in M^{(1)}$, the projective variety $\M_{m^1}$ is reducible, and $(\MM)_{m^1}$ and $(\MMM)_{m^1}$ are its irreducible components.

We conclude this preliminary part of the paper by observing that, locally, $M^{(k)}$ is the $k+1\St$ jet-extension of the trivial bundle $\mathbb{R}^2\times\mathbb{R}\to\mathbb{R}^2$, so that the reader more at ease with jet formalism may perform the substitution
 \begin{gather}\label{sostituzionePerIGettofili}
 M^{(k)}\longleftrightarrow J^{k+1}(2,1) ,\qquad k\geq 0 ,
\end{gather}
and easily adapt the main results to the new setting.

\section[Characteristics of $3\Rd$ order PDEs]{Characteristics of $\boldsymbol{3\Rd}$ order PDEs}\label{SecChar3ordPDEs}

Now we turn our attention to the hypersurfaces in $M^{(2)}$ that play the role of $3\Rd$ order PDEs in our analysis.

\begin{Definition} A characteristic covector $\alpha\in L_{m^2}^\ast{\setminus}\{0\}$, with $m^2\in\E_{m^1}$, is \emph{characteristic for~$\E$ in~$m^2$} if the corresponding rank-one line $\ell$ is tangent to $\E_{m^1}$ in $m^2$, in which case $\Span{\alpha}$ is a~\emph{characteristic $($direction$)$ for~$\E$ in~$m^2$} and the subspace~$H_\alpha$ is a \emph{characteristic hyperplane for~$\E$}.
\end{Definition}

Formula
 \eqref{eqGeometriaDelRaggio} says precisely that the rank-one line $\ell$ which corresponds to $\Span{\alpha^{\otimes 3}}$ via the isomorphism \eqref{isoFond} is precisely the (one-dimensional) tangent space to the prolongation $H_\alpha^{(1)}$. In this perspective, $\Span{\alpha}$ is a characteristic for $\E$ at $m^2$ if and only if $ H_\alpha^{(1)}$ is tangent to $\E_{m^1}$ at $m^2$ but, in general, $ H_\alpha^{(1)}$ does not need to touch $\E_{m^1}$ in any other point.
\begin{Definition}\label{defStrongChar}
 If $ H_\alpha^{(1)}$ is {entirely contained} in $\E_{m^1}$, then $\Span{\alpha}$ is called a \emph{strong characteristic $($direction$)$} and~$H_\alpha$ a~\emph{strongly characteristic line} (in~$m^2$).
\end{Definition}

\begin{Example}\label{exCharCoord}
Function \eqref{equation} determines, for $k=2$, the $3\Rd$ order PDE $\E=\{F=0\}$. Now equation \eqref{characteristic} can be correctly interpreted as follows: it is satisf\/ied if and only if the vector $\w=\w^i\xi_i$ given by \eqref{eq.char.dir} spans a characteristic line of $\E$. The very same equation \eqref{characteristic} tells also when the covector $\alpha=\w^2dx^1-\w^1dx^2$ spans a characteristic direction of $\E$. In the last perspective, equation \eqref{characteristic} is nothing but the right-hand side of \eqref{isoFondCoord} applied to $F$ and equated to zero.
\end{Example}

\subsection[Characteristics of a $3\Rd$ order PDE and relationship with its symbol]{Characteristics of a $\boldsymbol{3\Rd}$ order PDE and relationship with its symbol}\label{subsubDefPDE}

For any $m^2\in \E_{m^1}$ we def\/ine the \emph{vertical tangent space} to $\E$ at $m^2$ as the subspace
\begin{gather}\label{eqDefVertTangSpace}
V_{m^2}\E:=T_{m^2}\E_{m^1}\leq V_{m^2}M^{(2)},
\end{gather}
called the \emph{symbol} of $\E$ by many authors\footnote{Herewith we prefer to use the term ``symbol'' only for the function $F$ determining $\E$, and not for $\E$ itself. The two things, however, are the same, if one just regards $F$ as a the ($3\Rd$ order) dif\/ferential operator def\/ining the equation $\E$.}. Obviously, vertical tangent spaces can be naturally assembled into a linear bundle
\begin{gather*}
V\E:=\coprod_{m^2\in \E}V_{m^2}\E,
\end{gather*}
called \looseness=-1 the \emph{vertical bundle} of $\E$. Directly from \eqref{eqDefVertTangSpace} it follows the bundle embedding $V\E\subseteq \left.VM^{(2)}\right|_\E$,
where f\/ibres of the former are hyperplanes in the f\/ibres of the latter. Thanks to Lemma \ref{lemFund}, $V_{m^2}\E$ can also be regarded as a subspace of $S^3L_{m^2}^\ast$, and, being the identif\/ica\-tion~\eqref{isoFond} manifestly conformal, such an inclusion descends to the corresponding projective spaces,~i.e.,
\begin{gather}\label{InclusioneApparentementeStupida}
\P V_{m^2}\E \subseteq \P S^3L_{m^2}^\ast.
\end{gather}
The dual of canonical inclusion \eqref{InclusioneApparentementeStupida} reads
\begin{gather}\label{InclusioneApparentementeStupidaDUAL}
\Ann(V_{m^2}\E ) \in \P S^3L_{m^2},
\end{gather}
where $\Ann(V_{m^2}\E )$ is the (one-dimensional) subspace $\left(V_{m^2}M^{(2)}\right)^*$ made of covectors vanishing on the hyperplane $V_{m^2}\E $, viz.
\begin{gather}\label{eqComeSiDEfinisceIlSimbolo}
\Ann(V_{m^2}\E )\leq \big(V_{m^2}M^{(2)}\big)^*=(S^3L_{m^2}^\ast)^\ast=S^3L_{m^2}.
\end{gather}
From a global perspective, \eqref{InclusioneApparentementeStupidaDUAL} is nothing but the def\/inition of a section
\begin{gather}\label{secSymE}
\xymatrix{
\P S^3L \ar[r]& \E.\ar@/_/[l]_{\Ann(V\E )}
}
\end{gather}

\begin{Remark}[symbol of a function]\label{remSimboloSimboloso}
Let $\E=\{ F=0\}$, and use the same coordinates \eqref{coordELLE} of Remak \ref{remcoordELLE}. Then
 \begin{gather*}
\Ann(V_{m^2}\E ) = \Span{ \operatorname{Smbl}_{m^2}{F} },
\end{gather*}
where $\operatorname{Smbl}_{m^2}{F}$ is the \emph{symbol of $F$} at $m^2$, i.e.,
\begin{gather}\label{simbcordF}
\operatorname{Smbl}_{m^2}{F}=\left.\frac{\partial F}{\partial p_{ijk}}\right|_{m^2}\xi_i\xi_j\xi_k,
\end{gather}
where $\xi_i$ has been def\/ined in \eqref{defYi}.
Observe that $\Ann(V_{m^2}\E ) $ is independent of the choice of $F$ in the ideal determined by $\E$, so that $\Ann(V\E)$ can be replaced with $\operatorname{Smbl} {F}$.
\end{Remark}

In view of \eqref{eqComeSiDEfinisceIlSimbolo}, a direction $\Span{\alpha}\in \P L_{m^2}^\ast$ is a characteristic one for $\E$ at $m^2$ if and only if
 \begin{gather}\label{eqPairingSimpatico}
\big\langle \operatorname{Smbl}_{m^2}{F}, \alpha^{\otimes 3} \big\rangle=0,
\end{gather}
where $\langle \,\cdot\, ,\, \cdot\, \rangle$ is the canonical pairing on $S^3L_{m^2}$. Needless to say, \eqref{eqPairingSimpatico} is independent of the choice of $\alpha$ (resp., $F$) representing $\langle\alpha\rangle$ (resp., $\E$). Similarly, $\elle=\Span{\w}$, with $\w=\w^i\xi_i$, is a characteristic line if and only if
\begin{gather*}
 \sum_{i+j=3} (-1)^i \frac{\partial F }{\partial p_{\text{\scriptsize $\underbrace{1\cdots 1}_{i}$} \text{\scriptsize $\underbrace{2\cdots 2}_{j}$}}} \Bigg|_{ {m}^{2}}\big(\w^1\big)^i\big(\w^2\big)^j =0,
\end{gather*}
which eventually clarif\/ies formula \eqref{characteristic}.

\begin{Remark}[coordinates of $V\E$]
If $\E=\{F=0\}$ is a $3\Rd$ order PDE according to the def\/inition given at the end of Section~\ref{sec.definitions.contact.and.prol}, then $dF|_{\E_{m^1}}$ is nowhere zero, for all $m^1\in M^{(1)}$ (recall Section~\ref{sezChiarificatrice}). So,
 \begin{gather*}
V_{m^2}\E = \ker \left.\frac{\partial F}{\partial p_{ijk}}\right|_{m^2}dp_{ijk}
\end{gather*}
is an hyperplane in the 4D space $V_{m^2}M^{(2)}$.
\end{Remark}

\begin{Remark}[types of $3\Rd$ order PDEs]\label{RemTipiDiEquazioni}
Observe that~\eqref{simbcordF} is a $3\Rd$ order homogeneous polynomial in two variables, so that we can introduce its discriminiant
$\Delta=18abcd-4b^3d+b^2c^2-4ac^3-27a^2d^2$,
where $a,b,c,d\in C^\infty(M^{(2)})$ are def\/ined by
\begin{gather*}
a(m_2)=\left.\frac{\partial F}{\partial p_{111}}\right|_{m^2}\!,\qquad\!\!
b(m_2)=\left.\frac{\partial F}{\partial p_{112}}\right|_{m^2}\!,\qquad\!\!
c(m_2)=\left.\frac{\partial F}{\partial p_{122}}\right|_{m^2}\!,\qquad\!\!
d(m_2)=\left.\frac{\partial F}{\partial p_{222}}\right|_{m^2}\!.
\end{gather*}
Thus, any point $m_2\in\E$ can be of three dif\/ferent types, according to the signature of $\Delta$, namely,
\begin{itemize}\itemsep=0pt
\item $\Delta(m_2)>0$ $\Leftrightarrow$ \eqref{simbcordF} decomposes into three distinct linear factors;
\item $\Delta(m_2)=0$ $\Leftrightarrow$ \eqref{simbcordF} contains the square of a linear factor;
\item $\Delta(m_2)<0$ $\Leftrightarrow$ \eqref{simbcordF} contains an irreducible quadratic factor.
\end{itemize}
\end{Remark}

Equation $\E=\{F=0\}$ may be called ``fully parabolic'' in $m^2\in\E$ if \eqref{simbcordF} reduces to the cube of linear factor in the point $m^2$.

\subsection{The characteristic variety as a covering of the characteristic cone}
\label{sezionedovesidefinisceilconocaratteristico}
Let $\E=\{F=0\}$, and regard $L$ as a sub-bundle of the pull-back of $\C^1$ to the equation $\E$, bearing in mind the diagram
\begin{gather}\label{diagNonProjEQ}
\begin{split}&
\xymatrix{
 L \ar@{^(->}[r]^{}\ar[d]& \C^1\ar[d]\\
\E\ar[r] & M^{(1)}}\qquad \xymatrix{ L_{m^2} \ar@{^(->}[r]^{}\ar[d]& \C^1_{m^1}\ar[d]\\
**[l]\E\ni m^2\ar[r]^{\pi_{2,1}} &**[r] m^1\in M^{(1)}}
\end{split}
\end{gather}
Now we can come back to the direction $\w$ def\/ined by~\eqref{eq.char.dir} and notice that it lies in $L_{\overline{m}^2}$, and also provide an intrinsic way to check whether $\w$ is characteristic or not. Indeed, in view of the fundamental isomorphism~\eqref{isoFond},
one can regard the cube $\w^{\otimes 3}$ of $\w$ as a tangent vector to the f\/ibre $M^{(2)}_{\overline{m}^2}$, and equation~\eqref{characteristic} tells precisely when $\w^{\otimes 3}$ belongs to the sub-space $V_{\overline{m}^2}\E=T_{\overline{m}^2}\E_{\overline{m}^1}$ (see~\eqref{eqDefVertTangSpace}), up to a line-hyperplane duality (see also Remark~\ref{RemarcoFondamentale}). In other words, the set of characteristic lines
 \begin{gather*}
\Char{m^2}{\E}:=\{\elle\in \P L^\ast_{m^2}\,|\, \elle\textrm{ is a \emph{characteristic hyperplane} of }\E\}\subseteq \P L^\ast_{m^2}
\end{gather*}
is the projective sub-variety in $ \P L^\ast_{m^2}$ cut out by \eqref{characteristic}.

Since we plan to carry out an analysis of certain PDEs via their characteristics, it seems na\-tural to consider all characteristic lines at once, i.e., as a unique geometric object. Traditionally, one way to accomplish this is to take the disjoint union
\begin{gather*}
\Char{}{\E}:=\coprod_{m^2\in\E}\Char{m^2}{\E},
\end{gather*}
 known as the \emph{characteristic variety} of $\E$ (see, e.g., \cite{MR1083148}). As a bundle over $\E$, the family of the f\/ibres of $\Char{}{\E}$ coincides with the equation $\E$ itself. So, the bundle $\Char{}{\E}$, in spite of its importance for the study of a given equation $\E$, cannot be used to \emph{define} $\E$ as an object pertaining to~$M^{(1)}$. Still this can be arranged: it suf\/f\/ices ``to project everything one step down'', so to speak. The result is precisely the characteristic cone~$\V^\E$:
 \begin{gather}\label{eq.char.cone.2}
\underset{\textrm{Characteristic variety }\Char{}{\E}}{\boxed{\textrm{Projective sub-bundle of }\P L\to M^{(2)}}}\stackrel{(\pi_{2,1})_\ast}{\longrightarrow}\underset{\textrm{Characteristic cone }\M^{\E}}{\boxed{\textrm{Projective sub-bundle of }\P \C^1\to M^{(1)}}}\!\!\!
\end{gather}
It is easy to see that some characteristic lines, which are distinct entities in $\Char{}{\E}$, may collapse into $\M^{\E}$ (e.g., when the corresponding f\/ibres of the tautological bundle have non-zero intersection): hence, $\M^{\E}$ is not the most appropriate environment for the study of characteristics. On the other hand, $\M^{\E}$ is a bundle over $ M^{(1)}$, where there is no trace of the original equation $\E$: as such, it may ef\/fectively replace the equation itself and serves as a source of its invariants.

Now we can give a rigorous def\/inition of our main tool, the \emph{characteristic cone} of $\E$. Indeed (see Remark~\ref{RemarcoFondamentale}), the characteristic variety $\Char{}{\E}$ can be regarded as a sub-bundle of $\P L\longrightarrow \E$. In turn, this makes it possible to use the commutative diagram \eqref{diagNonProjEQ} to map $\Char{}{\E}$ to $\P\C^1$ and def\/ine the sub-bundle
\begin{gather*}
\M^\E:=\big\{ H\in \P\C^1\,|\, \exists \, m^2\in\E\colon {H}\subset L_{m^2}\textrm{ and}\\
\hphantom{\M^\E:=\big\{}{} H\textrm{ is a characteristic hyperplane of }\E\textrm{ at } m^2 \big\}\subseteq\P\C^1
\end{gather*}
 as the image of such a mapping.
\begin{Definition}\label{defCentrale}
Above def\/ined bundle $\M^\E\longrightarrow M^{(1)}$ is the \emph{characteristic cone} of $\E$.
\end{Definition}
By its def\/inition, $\M^\E$ f\/its into the commutative diagram
\begin{gather}\label{diagNonProjEQ-Char}
\begin{split}&
\xymatrix{ \Char{}{\E} \ar@{->>}[r]^{}\ar[d]& \M^\E\ar[d]\\
\E\ar[r] & M^{(1)}}
\end{split}
\end{gather}
revealing that, in a sense, $\M^\E$ is covered by $ \Char{}{\E}$.

It is worth observing that, from a local perspective,
\eqref{eq.char.cone.1}, \eqref{eq.char.cone.2}, and Def\/inition~\ref{defCentrale} all def\/ine the same object.
For instance, by Def\/inition \ref{defCentrale}, a line $\elle=\Span{\w}\in \P\C^1_{m^1}$
 belongs to $ \M^\E_{m^1}$ if and only if there is a point $m^2\in\E_{m^1}$ such that $\elle$ is a characteristic line for~$\E$ at the point~$m^2$. But Example~\ref{exCharCoord} shows that this is the case if and only if the generator $\w=(\w^1,\w^2)$ of $\elle$ satisf\/ies equation~\eqref{characteristic}.

 \subsection[The irreducible component $\MM^\E$ of the characteristic cone of a $3\Rd$ order PDE $\E$]{The irreducible component $\boldsymbol{\MM^\E}$ of the characteristic cone\\ of a $\boldsymbol{3\Rd}$ order PDE $\boldsymbol{\E}$}\label{subsubComponenteLineare}

Let $\E=\{F=0\}$ and identify $\Ann(V\E)$ with $\operatorname{Smbl} F$ as in Remark~\ref{remSimboloSimboloso}.
As a homogeneous cubic tensor on $L$ (see~\eqref{secSymE}), $\operatorname{Smbl} F$ possesses a linear factor, i.e., a section
\begin{gather*}
\xymatrix{
\P L\ar[r] &
\E\ar@/_/[l]_{\operatorname{Smbl}_I F}
}
\end{gather*}
such that
\begin{gather}\label{eqSpezzamentoSimbolo}
\operatorname{Smbl} F=\operatorname{Smbl}_I F\odot \operatorname{Smbl}_{II} F,
\end{gather}
 where $ \operatorname{Smbl}_{II} F$ is a section of $\P S^2L\rightarrow \E$.

 Plugging \eqref{eqSpezzamentoSimbolo} into \eqref{eqPairingSimpatico} shows that, in order to have
 \begin{gather*}
\big\langle \operatorname{Smbl}_{I,m^2}{ F}\odot \operatorname{Smbl}_{II,m^2} F, \alpha\odot\alpha^2 \big\rangle=\langle \operatorname{Smbl}_{I,m^2}{F},\alpha\rangle \big\langle \operatorname{Smbl}_{II,m^2}{F},\alpha^2\big\rangle=0
\end{gather*}
 it suf\/f\/ices that
 \begin{gather}\label{eqSimboletto}
\langle \operatorname{Smbl}_{I,m^2}{F},\alpha\rangle=0.
\end{gather}
In other words, if \eqref{eqSimboletto} is satisf\/ied, then $\Span{\alpha}\in \P L_{m^2}^\ast$ is a characteristic direction (incidentally revealing that $3\Rd$ order PDEs always possess a lot of them).

 Moreover, in view of the identif\/ication of characteristic directions with characteristic lines (Remark \ref{RemarcoFondamentale}), \eqref{eqSimboletto} shows also that the line $\operatorname{Smbl}_{I,m^2}{F}\in \P L_{m^2}$ is \emph{always} a characteristic line for $\E$, since there always is an $\alpha$ such that \eqref{eqSimboletto} is satisf\/ied.

It is clear now how diagram \eqref{diagNonProjEQ} can be made use of and how to def\/ine, much as we did in Section \ref{subsubRankUnoChar}, a subset $\MM^\E\subseteq\M^\E$ f\/itting into the commutative diagram
\begin{gather*}
\xymatrix{\operatorname{Smbl}_I F \ar[r] \ar[d]& \MM^\E\ar[d]\\
\E\ar[r] & M^{(1)}}
\end{gather*}
parallel to \eqref{diagNonProjEQ-Char}. This gives a solid background to the statement of Theorem~\ref{thEgr}, where~$\MM^\E$ was mentioned without exhaustive explanations.

\section{Proof of the main results}\label{secCentrale}

Now we are in position to prove the main results, Theorem \ref{LemmaNotturnoQuasiMattinieroPromossoATeorema}, concerned with the structure of $3\Rd$ order MAEs, Theorem~\ref{thEgr} and Theorem~\ref{thFormMag}, concerned with the structure of Goursat-type $3\Rd$ order MAEs. Since the former are special instances of the latter, we prefer to clarify their relationship beforehand.

 \subsection{Reconstruction of PDEs by means of their characteristics}\label{SubRicostruzione}
In spite of its early introduction \eqref{eq.char.cone.1}, the main object of our interest has been def\/ined again in a less direct way, which passes through the characteristic variety (see above Def\/inition~\ref{defCentrale}). The reason behind this choice is revealed by diagram \eqref{diagNonProjEQ-Char}: the characteristic variety $\Char{}{\E}$ plays the role of a~``minimal covering object'' for both the equation $\E$ itself and its characteristic cone $\V^\E$.

It is easy to guess\footnote{The structural relationship between the rank-one cone and the characteristic variety, which is the backbone of the present paper, made its appearance in many discussions concerning the characteristic variety, including \cite[Chapter~V]{MR1083148}, the late-1960's work by Guillemin on the characteristics of Pfaf\/f\/ian systems, and the more more recent work by Malgrange on the homological algebra of characteristic varieties.} that such a minimal covering must be def\/ined in terms of $\Omega$-isotropic f\/lags on~$M^{(1)}$, in the sense of Section~\ref{subMetaSim}, i.e., by means of the commutative diagram
 \begin{gather}\label{diagEVi}
 \begin{split}&
\xymatrix{& \Fl^{\textrm{iso}}_{2,1}(\C^1)\ar[dr]^{\p_1}\ar[dl]_{\p_2} & \\
M^{(2)}\ar@{.>>}[r] & M^{(1)} & \P\C^{(1)} \ar@{.>>}[l]}
\end{split}
\end{gather}
where $\Fl^{\textrm{iso}}_{2,1}(\C^1):=\{(m^2,\elle)\in M^{(2)}\times_{M^{(1)}}\P\C^1\,|\, L_{m^2}\supset\elle\}$
is the f\/lag bundle of $\Omega$-isotropic ele\-ments of $\C^1$.

The ``double f\/ibration transform'' associated to the diagram \eqref{diagEVi} allows to pass from a sub-bundle of $M^{(2)}$ (e.g., a $3\Rd$ order PDE $\E$) to a sub-bundle of $\P\C^1$ (e.g., its characteristic cone~$\M^\E$), and vice-versa. Lemma~\ref{lemFlagCharCone} indicates how to proceed in one direction.

\begin{Lemma}\label{lemFlagCharCone}
Let $\E$ be a PDE.
\begin{enumerate}\itemsep=0pt
\item[$1.$] Take the pre-image $\p_2^{-1}(\E)$ of $\E$ and select the points of tangency with the $\p_1$-fibres: the result is $\Char{}{\E}$.
\item[$2.$] Project $\Char{}{\E}$ onto $\P\C^1$ via $\p_1$: the result is~$\M^\E$.
\end{enumerate}
\end{Lemma}

\begin{proof}
 Item 1 follows from the fact that the $\p_1$-f\/ibres are the prolongations $H^{(1)}$ of hyperplanes $H\in\P\C^{(1)}$, and the tangency condition means that they are determined by characteristics (see~\eqref{eqGeometriaDelRaggio}). Item~2 is just a paraphrase of Def\/inition~\ref{defCentrale}.
\end{proof}
Coming back is easier. Indeed, given a sub-bundle $\M\subseteq\P\C^1$, the same def\/inition~\eqref{eqENu} of~$\E_\M$ given earlier can be recast as
 \begin{gather}\label{eqRicostruzione}
\E_\M=\p_2\big(\p_1^{-1}(\M)\big).
\end{gather}
In spite of the name ``double f\/ibration transform'', performing~\eqref{eqRicostruzione} f\/irst, and then applying Lemma~\ref{lemFlagCharCone} to the resulting equation~$\E_\M$, does not return, as a rule, the original sub-bundle~$\V$.
\begin{Corollary}\label{ColloarioCheServiraPrimaOPoi}
 Let $\M\subseteq\P\C^1$ be a sub-bundle. Then
\begin{gather}\label{eq.cuore}
\M\subseteq\M^{\E_\M}.
\end{gather}
\end{Corollary}
\begin{proof}
 Directly from Lemma \ref{lemFlagCharCone} and formula \eqref{eqRicostruzione}.
\end{proof}

We stress that, for any odd-order PDE $\E$, it holds the inclusion
\begin{gather}
\E\subseteq \E_{\M^\E}.\label{eqcuoricino}
\end{gather}
This paper begins to tackle the problem of determining those PDEs which are recoverable from their characteristics in the sense that \eqref{eqStar} is valid. Lemma \ref{LemmuccioSimpatico} below provides a simple evidence that, as a matter of fact, \emph{not all} PDEs are {recoverable} (in the sense of Def\/inition~\ref{defMaledetta}).

\begin{Lemma}\label{LemmuccioSimpatico}
 $\M$ is made of strongly characteristic lines for $\E_\M$.
\end{Lemma}
\begin{proof}
 Let $H\in \M$ be a line belonging to $\M$. Then $\p_2(\p_1^{-1}(H))$ is nothing but $H^{(1)}$ (see Def\/inition~\ref{defStrongChar}) and~\eqref{eqRicostruzione} tells precisely that $H^{(1)}$ is entirely contained in $\E_\M$. Hence, $H$ is a~strongly characteristic line for $\E_\M$ in any point of $\p_2(\p_1^{-1}(H))$.
\end{proof}

The meaning of Theorem \ref{thEgr} is that $3\Rd$ order MAEs of Goursat-type are, in a sense, those PDEs whose characteristic cone takes the simplest form, namely that of a linear projective sub-bundle, and, moreover, they are also recoverable from it.
Observe that \eqref{eqED} is a particular case of~\eqref{eqRicostruzione}, so that Goursat-type $3\Rd$ order MAEs constitute a remarkable example of equations determined by a projective sub-bundle of~$\P\C^1$. Then Lemma~\ref{LemmuccioSimpatico} reveals that all the lines lying in $\D$ are strongly characteristic lines for~$\E_\D$, in strict analogy with the classical case \cite{MR2985508}.

\subsection{Proof of Theorem \ref{LemmaNotturnoQuasiMattinieroPromossoATeorema}}

As we already underlined, we called ``$3\Rd$ order MAEs'' the equations \eqref{eqEOmega} def\/ined by means of a~2-form, disregarding the fact that, actually, Boillat introduced them by imposing the condition of complete exceptionality (which also allowed him to write down the local form~\eqref{eq.general.MAE.3}).

Fix now a point $m^1\in M^{(1)}$ and recall (see Lemma \ref{lemFund}) that the f\/ibre $M^{(2)}_{m^1}$ is an af\/f\/ine space modeled over $S^3L_{m^1}^\ast$. Let $\omega\in\Lambda^2\C^{1\,\ast}$ and regard the corresponding 2-form $\omega_{m^1}$ on $\C^1_{m^1}$ as a~linear map
$\Lambda^2\C^{1}_{m^1}\stackrel{\phi}{\longrightarrow} \R$.
The linear projective subspace
\begin{gather}\label{eqHypSec}
\mathcal{H}:=\P\ker\phi\subset \P \Lambda^2\C^{1}_{m^1}
\end{gather}
is the so-called \emph{hyperplane section} determined by $ \omega_{m^1}$. Now we show that $\mathcal{H}$ corresponds precisely to the f\/ibre $\E_{\omega, m^1}$ of the $3\Rd$ order MAE determined by $\omega$, via the Pl\"ucker embedding.

Indeed, $S^3L_{m^1}^\ast$ is embedded into the space $L_{m^1}^\ast\otimes S^2L_{m^1}^\ast$
 of \emph{all}, i.e., not necessarily Lagrangian, 2D horizontal subspaces, via the polarisation/Spencer operator. In turn, $L_{m^1}^\ast\otimes S^2L_{m^1}^\ast$ is an af\/f\/ine neighborhood of $L_{m^1}$ in the Grassmannian $\Gr(2,\C^1_{m^1})$, which is sent to $\P \Lambda^2\C^{1}_{m^1}$ by the Pl\"ucker embedding. Hence, the subspace $\E_{\omega, m^1}$ of $M^{(2)}_{m^1}$ can be regarded as a subspace of $\P \Lambda^2\C^{1}_{m^1}$, and \eqref{eqEOmega} tells precisely that such a subspace coincides with $\mathcal{H}$ def\/ined by \eqref{eqHypSec}.

We are now in position to generalize a result about classical MAEs \cite[Theorem 3.7]{MR2985508}, to the context of $3\Rd$ order MAEs.
\begin{Proposition}\label{proposizioneNotturna}
 A characteristic direction for $\E_\omega$ is also strongly characteristic.
\end{Proposition}
\begin{proof}
 Let $\elle\in\P\C^1_{m_1}$ be a characteristic direction for $\E_\omega$ at the point $m^2$. This means that $\elle\subset L_{m^2}$ and that the prolongation $\elle^{(1)}$ determined by $\elle$ is tangent to $\E_{\omega, m^1}$ at $m^2$. We need to prove that the whole $\elle^{(1)}$ is contained in~$\E_{\omega, m^1}$.

 To this end, recall that $\elle^{(1)}\subset M^{(2)}_{m^1}$ is a 1D af\/f\/ine subspace modeled over $S^3\Ann\elle$, passing through $L_{m^2}$ (see, e.g., the proof of Theorem~1 in~\cite{MorBach}). By the above arguments, $\elle^{(1)}$ can be embedded into $\P \Lambda^2\C^{1}_{m^1}$ as well. Now we can compare $\elle^{(1)}$ and $\mathcal{H}$: they are both linear, they pass through the same point~$L_{m^2}$, where they are also tangent each other. Hence, $\elle^{(1)}\subset \mathcal{H}$.
\end{proof}

Now Proposition \ref{proposizioneNotturna} allows us to prove Theorem \ref{LemmaNotturnoQuasiMattinieroPromossoATeorema}.

\begin{proof}[Proof of Theorem~\ref{LemmaNotturnoQuasiMattinieroPromossoATeorema}]
Inclusion \eqref{eq.cuore} is valid for \emph{any} $3\Rd$ order PDE. Conversely, if $m^2\in\E_{\M^\E}$, then there is a line $\elle\in \M^\E$, such that $L_{m^2}\supset\elle$. But
$\elle$ is a strong characteristic line for $\E$ thanks to
Proposition~\ref{proposizioneNotturna}, so that all Lagrangian planes passing through $\elle$ and, in particular, $L_{m^2}\equiv m^2$ itself, belong to $\E$.\end{proof}

\subsection{Proof of Theorem \ref{thEgr}}\label{subProofThEgr}

As it was outlined in Section \ref{subDescrMainRes}, Goursat-type MAEs are MAEs which correspond to decomposable forms, modulo a certain ideal. Before proving Theorem \ref{thEgr}, we make rigorous this statement. To this end, we shall need the submodule of \emph{contact $2$-forms}
\begin{gather}\label{defThetaGrande}
\Theta=\big\{\omega\in \Lambda^2 M^{(1)}\,|\,\omega|_{L_{m^2}}=0\ \forall\, m^2\in M^{(2)}\big\} \subset\Lambda^2 M^{(1)},
\end{gather}
and the corresponding projection
\begin{gather}\label{eqHor}
\Lambda^2 M^{(1)}\longrightarrow \frac{\Lambda^2 M^{(1)}}{\Theta}.
\end{gather}
Observe that the quotient bundle $ \frac{\Lambda^2 M^{(1)}}{\Theta}$ is canonically isomorphic to the rank-one bundle $\Lambda^2L^\ast$ over $M^{(1)}$. Hence, \eqref{eqHor} can be thought of as ($\Lambda^2L^\ast$)-valued.
\begin{Remark}\label{remarcoCheEraAddiritturaUnTheorema}
Two 2-forms have the same projection \eqref{eqHor} if and only if they dif\/fer by an element of $\Theta$.
\end{Remark}

\begin{Proposition}\label{proposizioneVeramenteSudata}
 For a $3D$ distribution $\D\subset\C^1$, it holds
 \begin{gather*}
\E_\D=\E_\omega \ \Leftrightarrow \ \omega=\rho_1\wedge\rho_2 \mod \Theta.
\end{gather*}
Moreover, the $3D$ sub-distribution $\D\subseteq\C^1$ is given by
\begin{gather}\label{eqDiComeKernel}
\D=\ker \rho_1|_{\C^1}\cap \ker \rho_2|_{\C^1}.
\end{gather}
\end{Proposition}
\begin{proof}
If $\Ann\D=\Span{\rho_1,\rho_2}\subset\C^{1\,\ast}$ is the annihilator of $\D$ in $\C^1$, then
$\E_\D$ can be written as~$\E_{\widetilde{\rho_1}\wedge\widetilde{\rho_2}}$, where $\widetilde{\rho_1},\widetilde{\rho_2}\in\Lambda^1M^{(1)}$ are extensions of $\rho_1$, $\rho_2$, respectively. The
 result follows from Remark~\ref{remarcoCheEraAddiritturaUnTheorema}.

Conversely, in light of \eqref{eqEOmega} and \eqref{defThetaGrande}, $\omega$ and $\rho_1\wedge\rho_2+\Theta$ give rise to the same equation, i.e.,
 \begin{gather}\label{eqEomegaErho12}
\E_\omega=\E_{\rho_1\wedge\rho_2}.
\end{gather}
It remains to be proved that the right-hand side of \eqref{eqEomegaErho12} is of the form $\E_\D$. To this end, it suf\/f\/ices to def\/ine $\D$ as in \eqref{eqDiComeKernel},
\begin{gather}\label{eqTroppoTildata}
\E_\D= \E_{\widetilde{\rho_1|_{\C^1}}\wedge\widetilde{\rho_2|_{\C^1}}},
\end{gather}
 and observe that the right-hand sides of \eqref{eqEomegaErho12} and
\eqref{eqTroppoTildata} coincide since
$
\widetilde{\rho_i|_{\C^1}}-\rho_i\in\Ann\C^1\subset\Theta$, $i=1,2$.
\end{proof}

Before starting the proof of Theorem \ref{thEgr}, we provide a meta-symplectic analog of the well-known formula $\Ann\D^\perp=\D\lrcorner\omega$.

\begin{Proposition}\label{propAnnullatoreEDirezioni}
 Let $\D_1\subseteq\C^1$ be a $3D$ sub-distribution with $\dim\D_1^v=2$ and $\E=\E_{\D_1}$ the corresponding quasi-linear $3\Rd$ order MAE. Then, if $\MMM^\E=\P\D_2\cup\P\D_3$,
 \begin{gather}
\Ann\D_i= {\frac{\D_1\lrcorner \Omega}{\elle_i}}\subseteq\C^{1\,\ast}, \qquad i=2,3, \label{eqAnnullDi}
\end{gather}
where $\elle_2$ and $\elle_3$ are the characteristic lines of
 $\D_1^v$.
\end{Proposition}

\begin{proof}
 Let
$\D_1=\Span{aD_1+bD_2, X_1,X_2}$,
with $X_1$ and $X_2$ as in \eqref{eqX1X2}. Then, the quasi-linear $3\Rd$ order MAE determined by $\D_1$ is
\begin{gather*}
\E_{\D_1}=ap_{111}+(b-a(k_1+k_2))p_{112}+(ak_1k_2-b(k_1+k_2))p_{122}+bk_1k_2p_{222}.
\end{gather*}
In order to obtain the right-hand side of \eqref{eqAnnullDi}, compute f\/irst
$
aD_1+bD_2\lrcorner \Omega = (adp_{11}+bdp_{12})\otimes \partial_{p_1}+(adp_{12}+bdp_{22})\otimes\partial_{p_2}
$,
and factor it by $\elle_i$,
\begin{gather}\label{eqLineaAggiuntiva}
aD_1+bD_2\lrcorner \Omega = (adp_{11}+(b-k_ia)dp_{12}-k_ib dp_{22})\otimes\partial_{p_1}\quad \mod \elle_i.
\end{gather}
Combining \eqref{eqLineaAggiuntiva} above with \eqref{eqLineaFinale}, yields
\begin{gather}\label{eqD1OmegaModElli}
{\frac{\D_1\lrcorner \Omega}{\elle_i}}\simeq\Span{k_jdx^1+dx^2, adp_{11}+(b-k_ia)dp_{12}-k_ib dp_{22}}.
\end{gather}
Finally, the wedge product of the two 1-forms spanning the module~\eqref{eqD1OmegaModElli} above, is the 2-form
\begin{gather*}
\omega= ak_jdp_{11}\wedge dx^1+(b-k_ia)k_jdp_{12}\wedge dx^1-bk_ik_jdp_{22}\wedge dx^1 \\
\hphantom{\omega=}{} + adp_{11}\wedge dx^2+(b-k_ia)dp_{12}\wedge dx^2-k_ibdp_{22}\wedge dx^2,
\end{gather*}
and direct computations show that
$
\E_\omega=\E_{\D_1}
$.
The result follows from Proposition \ref{proposizioneVeramenteSudata}.
\end{proof}

Now we turn back to Theorem \ref{thEgr}, and deal separately with its two implications.

\subsubsection{Proof of the suf\/f\/icient part of Theorem \ref{thEgr}}

The suf\/f\/icient part of Theorem \ref{thEgr} will be proved through Lemma \ref{LemQLsonoMonge} below.

\begin{Lemma}\label{LemQLsonoMonge}
 Let $\E=\{F=0\}$, where $F$ is given either by \eqref{eqGoursat},
 or
 \begin{gather}\label{eqGenQL}
F=ap_{111}+bp_{112}+cp_{122}+dp_{222}+e,
\end{gather}
where $a,b,c,d,e\in C^\infty(M^{(1)})$.
Then there exists a $3D$ sub-distribution $\D\subseteq\C^{(1)}$ such that $\M^\E=\MM^\E\cup\MMM^\E$, with $\MM^\E=\P\D$. Moreover, $\E=\E_\D$, i.e., according to~\eqref{eqED}, it is a~$3\Rd$ order MAE of Goursat-type.
\end{Lemma}

We recall that Lemma \ref{LemQLsonoMonge} is to be interpreted in f\/ibre-wise perspective (see Section~\ref{sezChiarificatrice}), i.e., as a collection of statements, each of which corresponds to the f\/ixation of a point $m^1\in M^{(1)}$. But, since the formula~\eqref{eqGenQL} must correspond to an af\/f\/ine hyperplane in the space of the $p_{ijk}$'s over $m^1$, then the vector $(a(m^1),b(m^1),c(m^1),d(m^1))\in\R^4 $ must be non-zero. This means that for any $m^1\in M^{(1)}$, there is always one of the functions $a,b,c,d$, which is non-zero at $m^1$.

\begin{proof}[Proof of Lemma \ref{LemQLsonoMonge}]
Let us begin with the non-linear case.

In order to verify that $\E=\E_\D$, it suf\/f\/ices to put
\begin{gather}
\D:= \big\langle D_1 + f_{111}\partial_{p_{11}} + f_{112}\partial_{p_{12}} + f_{122}\partial_{p_{22}}, D_2 + f_{211}\partial_{p_{11}} + f_{212}\partial_{p_{12}} + f_{222}\partial_{p_{22}},\nonumber \\
\hphantom{\D:=\big\langle}{} R\partial_{p_{11}}+S\partial_{p_{12}}+T\partial_{p_{22}}\big\rangle,\label{eq.D.Goursat.nonlin}
\end{gather}
where $(R,S,T)=\boldsymbol{A}$ is the same appearing in \eqref{eqGoursat}, and
observe that (see also \eqref{eqCoordEll})
 \begin{gather*}
L_{m^2}=\Span{ D_1+ p_{111}\partial_{p_{11}}+p_{112}\partial_{p_{12}}+p_{122}\partial_{p_{22}}, D_2 + p_{112}\partial_{p_{11}}+p_{122}\partial_{p_{12}}+p_{222}\partial_{p_{22}} }
\end{gather*}
belongs to $\E_{\D}$ (see \eqref{eqED}) if and only if the $5\times 5$ determinant
\begin{gather*}
\det\left|\begin{matrix}
1 & 0 & p_{111} & p_{112} & p_{122}
\\
0 & 1 & p_{112} & p_{122} & p_{222}
\\
1 & 0 & f_{111} & f_{112} & f_{122}
\\
0 & 1 & f_{211} & f_{212} & f_{222}
\\
0 & 0 & R & S & T
\end{matrix}\right|
=
\det \left|
\begin{matrix}
0 & 0 & p_{111}-f_{111} & p_{112} -f_{112}& p_{122} -f_{122}
\\
0 & 0 & p_{112}-f_{211} & p_{122}-f_{212} & p_{222}-f_{222}
\\
1 & 0 & f_{111} & f_{112} & f_{122}
\\
0 & 1 & f_{211} & f_{212} & f_{222}
\\
0 & 0 & R & S & T
\end{matrix}\right|=-F
\end{gather*}
is zero.

We now prove that $\MM^\E=\P\D$ in the case the coef\/f\/icient of $p_{111}$ of equation~\eqref{eqGoursat} is not zero. The other cases are formally analog. So, let us solve equation~\eqref{eqGoursat} with respect to~$p_{111}$ and take the remaining coordinates $p_{112}$, $p_{122}$, $p_{222}$ as local coordinates on~$\E$. Thus~$\operatorname{Smbl} F$ (see~\eqref{simbcordF}) is, up to a factor, equal to
\begin{gather}
\Big(\big(S(p_{222}-f_{222})-T(p_{122}-f_{212})\big)\xi_1+\big(T(p_{112}-f_{112})-S(p_{122}-f_{122})\big)\xi_2\Big)\nonumber \\
\qquad{}\times \Big(\big(T(p_{122}-f_{212})-S(p_{222}-f_{222})\big)\xi_1^2
+\big(R(p_{222}-f_{222})-T(p_{112}-f_{211})\big)\xi_1\xi_2\nonumber\\
\qquad{}+\big(S(p_{112}-f_{211})-R(p_{122}-f_{212})\big)\xi_2^2\Big),\label{eqSimboloPazzesco}
\end{gather}
so that (see also Section~\ref{subsubComponenteLineare})
\begin{gather}
 \MM^\E =\Big\{\Big(\big(S(p_{222}-f_{222})-T(p_{122}-f_{212})\big)D_1
 \nonumber\\
\hphantom{\MM^\E =}{} + \big(T(p_{112}-f_{112})-S(p_{122}-f_{122})\big)D_2\Big)_{m^2} \,|\, m^2\in\E\Big\}.\label{eq.V1.Goursat.nonlin}
\end{gather}
A direct computation shows that
\begin{gather*}
 \MM^\E=\big\{\big(T(p_{122}-f_{212})-S(p_{222}-f_{222})\big)(D_1 + f_{111}\partial_{p_{11}} + f_{112}\partial_{p_{12}} + f_{122}\partial_{p_{22}})
\\
\hphantom{\MM^\E=}{} +
\big(S(p_{122}-f_{122})-T(p_{112}-f_{112})\big)(D_2 + f_{211}\partial_{p_{11}} + f_{212}\partial_{p_{12}} + f_{222}\partial_{p_{22}})
\\
\hphantom{\MM^\E=}{}
+ (-p_{112}p_{222}+p_{112}f_{222}+f_{112}p_{222}-f_{112}f_{222}
+p_{122}^2-p_{122}f_{212}-f_{122}p_{122}+f_{122}f_{212})\! \\
\hphantom{\MM^\E=}{}
\times (R\partial_{p_{11}}+S\partial_{p_{12}}+T\partial_{p_{22}})\,|\,
p_{112},p_{122},p_{222}\in\mathbb{R}
\big\}
\end{gather*}
so that $ \MM^\E$ turns out to be the 3D linear space~\eqref{eq.D.Goursat.nonlin}.

Let us now pass to the quasi-linear case.

One of the coef\/f\/icients of the third derivatives of \eqref{eqGenQL} must be non-zero. Assume that $a\neq 0$; the remaining cases can be treated similarly. Again, we solve equation~\eqref{eqGenQL} with respect to~$p_{111}$, so that $p_{112}$, $p_{122}$, $p_{222}$ become local coordinates on it, and
\begin{gather} \label{eq-symb-quasi-linear}
\operatorname{Smbl} F=a\xi_1^3+b\xi_1^2\xi_2+c\xi_1\xi_2^2+d\xi_2^3=(k\xi_1+h\xi_2)q(\xi_1,\xi_2),
\end{gather}
where $q(\xi_1,\xi_2)$ is a homogeneous quadratic function in $\xi_1$ and $\xi_2$. Thus, we have that (see again Section~\ref{subsubComponenteLineare}).
\begin{gather} \label{eq-MME-Goursat-linear}
\MM^\E=\Span{ (kD_1+hD_2)_{m^2} \,|\, m^2\in\E}.
\end{gather}
Let us assume $d=0$.
In order to verify that $\E=\E_\D$, it is enough to put
\begin{gather}\label{eq-D-Goursat-lin-d-zero}
\D:= \Span{ D_1-\frac{e}{a}\partial_{p_{11}}, -\frac{b}{a}\partial_{p_{11}}+\partial_{p_{12}}, -\frac{c}{a}\partial_{p_{11}}+\partial_{p_{22}}}.
\end{gather}
To prove that $\MM^\E=\P\D$, observe that the symbol \eqref{eq-symb-quasi-linear} contains the linear factor $\xi_1$, and
\begin{align*}
\MM^\E&=\langle {D_1}|_{m^2} \rangle_{m^2\in\E } \\
&=\left\{D_1-\frac{e}{a}\partial_{p_{11}}+p_{112}\left(-\frac{b}{a}\partial_{p_{11}}+\partial_{p_{12}}\right) + p_{122}\left(-\frac{c}{a}\partial_{p_{11}}+\partial_{p_{22}}\right)\,|\, p_{112},p_{122}\in\mathbb{R}\right\}
\end{align*}
is the projectivization of~\eqref{eq-D-Goursat-lin-d-zero}.

Finally assume $d\neq 0$.

Observe that the lines in $\MM^\E$, which, in view of \eqref{eq-MME-Goursat-linear}, are generated by
\begin{gather*}
 kD_1+hD_2-\frac{ke}{a}\partial_{p_{11}}+
p_{112}\left(\left(-\frac{kb}{a}+h\right)\partial_{p_{11}}+k\partial_{p_{12}}\right) \\
\qquad {} + p_{122}\left(-\frac{kc}{a}\partial_{p_{11}}+h\partial_{p_{12}}+k\partial_{p_{22}}\right)
+p_{222}\left(-\frac{kd}{a}\partial_{p_{11}}+h\partial_{p_{22}}\right),
\end{gather*}
$p_{112},p_{122},p_{222}\in\mathbb{R}$,
f\/ill a 3D linear space since $\frac{\xi_1}{\xi_2}=-\frac{h}{k}$ is a solution to \eqref{eq-symb-quasi-linear}, i.e.,
\begin{gather}\label{eq.abcde}
-a\frac{h^3}{k^3}+b\frac{h^2}{k^2}-c\frac{h}{k}+d=0.
\end{gather}
So, if we put
\begin{gather*}
\D:= \Span{ D_1 + \frac{h}{k}D_2 - \frac{e}{a}\partial_{p_{11}}, \left(-\frac{b}{a}+\frac{h}{k}\right)\partial_{p_{11}}+\partial_{p_{12}}, \frac{k}{h}\frac{d}{a}\partial_{p_{11}}-\partial_{p_{22}}}
\end{gather*}
with $h,k\neq 0$ such that \eqref{eq.abcde} is satisf\/ied, we obtain $\E=\E_\D$ and $\MM^\E=\P\D$.
\end{proof}

\subsubsection{Proof of the necessary part of Theorem \ref{thEgr}}\label{subDimThEgrNec}

Being $\E_{m^1}$ a closed submanifold of codimension~1, in the neighborhood of any point $m^2\in \E_{m^1}$, we can always present $\E$ in the form $\E=\{F=0\}$, with $F:= p_{ijk}-G$, with $G$ not depending on~$p_{ijk}$, for some $(i,j,k)$ which we shall assume equal to~$(1,1,1)$, since the other cases are formally analog. Then the symbol of $F$ at $m^1$ is
\begin{gather}\label{eqSymb1}
\operatorname{Smbl}_{m^1}F=\xi_1^3-G_{p_{112}}\xi_1^2\xi_2-G_{p_{122}}\xi_1\xi_2^2-G_{p_{222}}\xi_2^3.
\end{gather}
The right-hand side of \eqref{eqSymb1} is a $3\Rd$ order homogeneous polynomial with unit leading coef\/f\/icient: hence, there exist unique~$\beta$,~$A$,~$B$ such that
$\operatorname{Smbl}_{m^1}F=(\xi_1+\beta\xi_2)(\xi_1^2+A\xi_1\xi_2+B\xi_2^2)$.
Following the general procedure (see also formula~\eqref{eq.char.cone.2} and Def\/inition~\ref{defCentrale}), to construct the characteristic cone~$\M^\E$, one easily sees that $\M^\E$ contains
 the following 3-parametric family of lines:
 \begin{gather*}
\MM^\E:=\big\{ D_1+\beta D_2+ (G+\beta p_{112})\partial_{p_{11}}+ (p_{112}+\beta p_{122})\partial_{p_{12}}\\
\hphantom{\MM^\E:=\big\{}{} + (p_{122}+\beta p_{222})\partial_{p_{22}} \,|\, p_{112},p_{122},p_{222}\in \R \big\}.
\end{gather*}
Observe that, with the same notation as~\eqref{eq.char.dir}, the direction $v=(v^1,v^2)$ belongs to~$\MM^\E$ if and only if there exists $m^2\in\E_{m^1}$ such $v^1+\beta (m^2) v^2=0$ (see also Remark~\ref{RemarcoFondamentale}).

Suppose now that there exists a 3D sub-distribution $\D\subset\C ^1$ such that
\begin{gather}\label{eqCondizTh1}
\MM^\E=\P \D.
\end{gather}
Since $\D$ has codimension 2 in $\C ^1$, condition \eqref{eqCondizTh1} can be dualized as follows: there are two independent forms $\rho_1,\rho_2\in\C ^{1\, \ast}$,
\begin{gather*}
\rho_1 = k_1dx^1+k_2dx^2+k^{11}dp_{11}+k^{12}dp_{12}+k^{22}dp_{22},
\\
\rho_2 = h_1dx^1+h_2dx^2+h^{11}dp_{11}+h^{12}dp_{12}+h^{22}dp_{22},
\end{gather*}
such that $\Ann\D=\langle\rho_1,\rho_2\rangle$, i.e., $\D=\ker\rho_1\cap\ker\rho_2$. In other words, \eqref{eqCondizTh1} holds true if and only if $\rho_1(v)=\rho_2(v)=0$
for all $v\in\MM^\E$, i.e.,
\begin{gather}
k_1+k^{12}p_{112}+k^{22}p_{122}+\big(k_2+k^{11}p_{112}+k^{12}p_{122}+k^{22}p_{222}\big)\beta+k^{11}G = 0,\label{eqSystema1}\\
h_1+h^{12}p_{112}+h^{22}p_{122}+\big(h_2+h^{11}p_{112}+h^{12}p_{122}+h^{22}p_{222}\big)\beta+h^{11}G = 0,\label{eqSystema2}
\end{gather}
identically in $p_{112}$, $p_{122}$, $p_{222}$,
where the $2\times 5$ matrix{\samepage
 \begin{gather}\label{matricechenonpiaceagianni}
\left|\begin{matrix}k_1 & k_2 & k^{11} & k^{12} & k^{22} \\h_1 & h_2 & h^{11} & h^{12} & h^{22}\end{matrix}\right|
\end{gather}
of functions on $M^{(1)}$ has rank two everywhere.}

Now \eqref{eqSystema1}, \eqref{eqSystema2} must be regarded as linear system of two equations in the unknowns~$\beta$ and~$G$. Its discriminant is easily computed,
$\Delta:=(h^{11}k^{22}-k^{11}h^{22})p_{222}+(h^{11}k^{12}-k^{11}h^{12})p_{122}+h^{11}k_2-k^{11}h_2$,
and $\Delta$ is a polynomial in $p_{122}$, $p_{222}$.

Suppose $\Delta(m^2)=0$, and interpret \eqref{eqSystema1}, \eqref{eqSystema2} as a linear system in the variables $\beta$ and $G$. The condition $\Delta(m^2)=0$, together with the compatibility conditions of such a system, implies that the matrix~\eqref{matricechenonpiaceagianni} is of rank $1$ in $m^1$. This contradicts the hypothesis that $\dim\D_{m^1}=3$, so that~$\Delta$ must be everywhere non-zero.

Then, in the points with $\Delta(p_{122},p_{222})\neq 0$ one can f\/ind $\Delta G$ as a polynomial expression of~$p_{112}$, $p_{122}$, $p_{222}$, whose coef\/f\/icients turn out to be minors of the matrix \eqref{matricechenonpiaceagianni}. In particular, $p_{111}-G$ can be singled out from the so-obtained expression, namely
\begin{gather}
-\Delta(p_{111}-G) =-h_1k_2+k_1h_2-\big(h^{11}k_2-k^{11}h_2\big)p_{111} \nonumber\\
\hphantom{-\Delta(p_{111}-G) =}{} -\big(k^{11}h_1-h^{11}k_1+k_2h^{12}-k^{12}h_2\big)p_{112}\nonumber\\
\hphantom{-\Delta(p_{111}-G) =}{}-\big(k_2h^{22}+k^{12}h_1-k^{22}h_2-k_1h^{12}\big)p_{122}-\big(k^{22}h_1-k_1h^{22}\big)p_{222} \nonumber\\
\hphantom{-\Delta(p_{111}-G) =}{}+\big(k^{12}h^{22}-h^{12}k^{22}\big)\big(p_{112}p_{222}-p_{122}^2\big) \nonumber\\
\hphantom{-\Delta(p_{111}-G) =}{} + \big(k^{22}h^{11}-k^{11}h^{22}\big)(-p_{111}p_{222}+p_{112}p_{122}) \nonumber\\
\hphantom{-\Delta(p_{111}-G) =}{}+\big(k^{11}h^{12}-h^{11}k^{12}\big) \big(p_{111}p_{122}-p_{112}^2\big).\label{eqDeltaPi111MenoGi}
\end{gather}
We need to show that \eqref{eqDeltaPi111MenoGi} is satisf\/ied if and only if there is a nowhere zero factor $\lambda$ such that
 \begin{gather}\label{eqFp111G}
F=\lambda (p_{111}-G),
\end{gather}
where $F$ is given by \eqref{eqGoursat}, i.e., equation~\eqref{eqFp111G} needs to be solved with respect to~$f_{ijk}$,~$R$, $S$,~$T$. By equating the coef\/f\/icients of~$p_{111}$, one obtains
\begin{gather}\label{eqLambda}
\lambda= T(p_{122}-f_{212})-S(p_{222}-f_{222}),
\end{gather}
and replacing~\eqref{eqLambda} into~\eqref{eqFp111G} yields
\begin{gather}\label{eqFTSG}
F-( T(p_{122}-f_{212})-S(p_{222}-f_{222}))(p_{111}-G)=0.
\end{gather}
Direct computations show that the left-hand side of
 \eqref{eqFTSG} is a rational function whose numerator is a $2\Nd$ order polynomial in $p_{112}$, $p_{122}$, $p_{222}$. From the vanishing of this polynomial, one can express $f_{ijk}$, $R$, $S$, $T $ as functions of the entries of the matrix \eqref{matricechenonpiaceagianni}.

\subsection{Proof of Theorem \ref{thFormMag}}
It will be accomplished in steps. As a preparatory result, we show that a 3D sub-distribution of $\C^1$ with 2D vertical part determines a quasi-linear $3\Rd$ order MAE (Lemma \ref{LemEd1QLMonge} below). The converse statement, i.e., that the characteristic cone of a generic quasi-linear $3\Rd$ order PDE possesses a linear sheet determined by a 3D sub-distribution $\D\subseteq \C^1$ with 2D vertical part, has been proved by Lemma \ref{LemQLsonoMonge} above. Then we pass to generic 3D sub-distributions of $\C^1$, and derive the expression of the corresponding equation $\E_\D$ (Lemma \ref{Lemma-ex-Considerazioni-Forma-Generale}). As a consequence of these results (Corollary \ref{CorollarioCheSeFunzionaSarebbeCarino}), we prove the initial statement of Theorem \ref{thFormMag}.

The next Sections \ref{SubStat1}, \ref{SubStat2}, \ref{sec.3.3} are devoted to the specif\/ic proofs of the items~1,~2,~3 of Theorem~\ref{thFormMag}, respectively.

Besides the proof of Theorem~\ref{thFormMag} itself, a few interesting byproducts will be pointed out. For instance, Corollary~\ref{CorollarioCheSeFunzionaSarebbeCarino} means that for quasi-linear~$3\Rd$ order MAEs, all characteristic lines are strong, generalizing an analogous result for classical MAEs (see \cite[Theorem~3.7]{MR2985508}). Also Corollary~\ref{corRicettaMagica} is a non-trivial and unexpected generalization of a phenomenon f\/irstly observed in the classical case (see~\cite[Theorem~1.1]{MR2985508}).

\begin{Lemma}\label{LemEd1QLMonge}
 Let $\D_1=\Span{h_1,\D_1^v}$, with
\begin{gather*}
h_1 = aD_1+bD_2+f_{11}\partial_{p_{11}}+f_{12}\partial_{p_{12}}+f_{22}\partial_{p_{22}},
\\
\D_1^v = \Span{R_i\partial_{p_{11}}+S_i\partial_{p_{12}}+T_i\partial_{p_{22}}\,|\, i=1,2}.
\end{gather*}
Then $\E_{\D_1}$ is of the form \eqref{eq.general.MAE.3} with
\begin{gather}
\boldsymbol{A} = 0,\label{condA}\\
C = \det \left|\begin{matrix}f_{11} & f_{12} & f_{22} \\R_1 & S_1 & T_1 \\R_2 & S_2 & T_2\end{matrix}\right|,\label{condC}\\
\boldsymbol{B} = \left(-a \det\left|\begin{matrix}S_1 & T_1 \\S_2 & T_2\end{matrix} \right| , -b \det\left|\begin{matrix}S_1 & T_1 \\S_2 & T_2\end{matrix}\right| +a \det\left|\begin{matrix}R_1 & T_1 \\R_2 & T_2\end{matrix}\right| ,\right.\nonumber\\
\left. \hphantom{\boldsymbol{B} =}{} \quad
 b \det\left|\begin{matrix}R_1 & T_1 \\R_2 & T_2\end{matrix}\right| - a \det\left|\begin{matrix}R_1 & S_1 \\R_2 & S_2\end{matrix}\right|, -b \det\left|\begin{matrix}R_1 & S_1 \\R_2 & S_2\end{matrix}\right| \right).\label{condB}
\end{gather}
\end{Lemma}

\begin{proof}
 Just observe that
\begin{gather*}
L_{m^2}=\Span{ D_1+ p_{111}\partial_{p_{11}}+p_{112}\partial_{p_{12}}+p_{122}\partial_{p_{22}}, D_2 + p_{112}\partial_{p_{11}}+p_{122}\partial_{p_{12}}+p_{222}\partial_{p_{22}} }
\end{gather*}
belongs to $\E_{\D_1}$ if and only if the determinant of the $5\times 5$ matrix
\begin{gather}\label{bigDET}
\left|\begin{matrix}a & b & f_{11} & f_{12} & f_{22} \\0 & 0 & R_1 & S_1 & T_1 \\0 & 0 & R_2 & S_2 & T_2 \\1 & 0 & p_{111} & p_{112} & p_{122} \\0 & 1 & p_{112} & p_{122} & p_{222}\end{matrix}\right|
\end{gather}
is zero (see also \eqref{eqCoordEll} and \eqref{eqED}). Standard matrix manipulations reveal that the determinant of~\eqref{bigDET} coincides in turn with \eqref{eq.general.MAE.3} under conditions \eqref{condA}, \eqref{condC} and \eqref{condB}.
\end{proof}

\begin{Remark}\label{remOttimale}
Since $\dim\Gr(3,5)=6$, it takes 6 parameters from $C^\infty(M^{(1)})$ to identify a 3D sub-distribution $\D\subset\C^1$ and, hence, an equation of the form~$\E_\D$.
\end{Remark}

Above Remark \ref{remOttimale} shows that the description \eqref{eqGoursat} of the $\E_\D$'s has some redundancies, and Lemma~\ref{lemma.optimal}) below
is a way to ref\/ine it.

\begin{Lemma}
\label{Lemma-ex-Considerazioni-Forma-Generale}
The local form of a generic $\E_\D$ with $\D\in\Gr(3, \C^1)$, in the case $\dim\D^v=1$, is given by~\eqref{eqGoursat}.
\end{Lemma}
\begin{proof}
A distribution $\D\in\Gr(3, \C^1)$ with $\dim\D^v=1$ is locally given by
\begin{gather*}
\D=\big\langle D_1+f_{111}\partial_{p_{11}} +f_{112}\partial_{p_{12}} +f_{122}\partial_{p_{22}} , D_2+f_{211}\partial_{p_{11}} +f_{212}\partial_{p_{12}} +f_{222}\partial_{p_{22}}, \\
\hphantom{\D=\big\langle}{} R\partial_{p_{11}} +S\partial_{p_{12}} +T\partial_{p_{22}} \big\rangle,
\end{gather*}
where $f_{ijk}$, $R$, $S$, $T$ are $C^\infty$ functions (def\/ined in some neighborhood of~$M^{(1)}$). Now it is clear that a Lagrangian plane
\begin{gather*}
\langle D_1+p_{111}\partial_{p_{11}} +p_{112}\partial_{p_{12}} +p_{122}\partial_{p_{22}} , D_2+p_{112}\partial_{p_{11}} +p_{122}\partial_{p_{12}} +p_{222}\partial_{p_{22}} \rangle
\end{gather*}
non-trivially intersects $\D$ if\/f
\begin{gather*}
\det
\left(
\begin{matrix}
1 & 0 & p_{111} & p_{112} & p_{122}
\\
0 & 1 & p_{112} & p_{122} & p_{222}
\\
1 & 0 & f_{111} & f_{112} & f_{122}
\\
0 & 1 & f_{211} & f_{212} & f_{222}
\\
0 & 0 & R & S & T
\end{matrix}
\right)
=
\left(
\begin{matrix}
p_{111}-f_{111} & p_{112}-f_{112} & p_{122}-f_{122}
\\
p_{112}-f_{211} & p_{122}-f_{212} & p_{222}-f_{222}
\\
R & S & T
\end{matrix}
\right)=0
\end{gather*}
that is equal to~\eqref{eqGoursat} with $\A=(R,S,T)$.
\end{proof}

Now we discuss the possibility of reducing the number of redundant parameters in~\eqref{eqGoursat} (see Remark~\ref{remOttimale}). In turn, this is linked to the natural question whether there exists a Lagrangian horizontal part $\mathcal{H}$ of $\mathcal{D}\in\Gr(3,\C^1)$ in the case that~$\dim\mathcal{D}^v=1$, i.e., the existence of a splitting $\mathcal{D}=\mathcal{H}\oplus\mathcal{D}^v$, with Lagrangian~$\mathcal{H}$.

\begin{Lemma}\label{lemma.appoggio.gianni}
Let $\mathcal{D}\subset\mathcal{C}^1$ be a $3D$ distribution such that $\dim\mathcal{D}^v=1$. If there exists a~Lag\-rangian horizontal part $\mathcal{H}$ of $\mathcal{D}$, then $\mathcal{E}_{\mathcal D}$ can be put in the form $\E=\{F=0\}$, where $F$ is given by~\eqref{eqGoursat} with $f_{112}=f_{211}$ and $f_{122}=f_{212}$.
\end{Lemma}

\begin{proof}
It follows easily taking into account the computations of Lemma \ref{Lemma-ex-Considerazioni-Forma-Generale}.
\end{proof}

\begin{Lemma}\label{lemma.optimal}
Let $\mathcal{D}\subset\mathcal{C}^1$ be a $3D$ distribution such that $\dim\mathcal{D}^v=1$. Let $m^1\in M^{(1)}$. If $\rank(\mathcal{D}^v_{m^1})\neq 1$ $($see Definition~{\rm \ref{defRankOneLine}} concerning rank-one lines$)$, then $\mathcal{E}_{\mathcal D}$ can be put, in a~neighborhood of~$m^1$, in the form $\E=\{F=0\}$, where $F$ is given by~\eqref{eqGoursat} with $f_{112}=f_{211}$ and $f_{122}=f_{212}$.
\end{Lemma}

\begin{proof}
Let the vertical part $\mathcal{D}^v$ of $\mathcal{D}$ be spanned by
\begin{gather*}
\mathcal{V}=R\partial_{p_{11}}+S\partial_{p_{12}}+T\partial_{p_{22}}
\end{gather*}
where $R,S,T\in C^\infty(M^{(1)})$. Let
\begin{gather*}
\mathcal{D}=
\langle D_1 + b_{111}\partial_{p_{11}} + b_{112}\partial_{p_{12}} + b_{122}\partial_{p_{22}},\,
D_2 + b_{211}\partial_{p_{11}} + b_{212}\partial_{p_{12}} + b_{222}\partial_{p_{22}},
\mathcal{V} \rangle,
\end{gather*}
with $b_{ijk}\in C^\infty(M^{(1)})$. For any $\alpha,\beta\in C^\infty(M^{(1)})$ we have that
\begin{gather}
\D =\big\langle D_1 + b_{111}\partial_{p_{11}} + b_{112}\partial_{p_{12}} + b_{122}\partial_{p_{22}}
+ \alpha \mathcal{V},\nonumber\\
\hphantom{\D =\big\langle}{}
D_2 + b_{211}\partial_{p_{11}} + b_{212}\partial_{p_{12}} + b_{222}\partial_{p_{22}}
+ \beta \mathcal{V} ,
\mathcal{V} \rangle.\label{eq.D.new}
\end{gather}
If $\rank(\mathcal{V}_{m^1})\neq 1$, then $\alpha$ and $\beta$ can be chosen in such a way that, in a neighborhood of~$m^1$,
\begin{gather}\label{eq.sys.A.B}
b_{112}+\alpha S=b_{211}+\beta R , \qquad b_{122}+\alpha T=b_{212}+\beta S.
\end{gather}
In fact system~\eqref{eq.sys.A.B} is always compatible for any
$R$, $S$, $T$, $b_{112}$, $b_{211}$, $b_{122}$, $b_{212}$ since \mbox{$\rank(\mathcal{V}_{m^1}){\neq} 1$} if\/f $RT-S^2\neq 0$ at~$m^1$. To conclude this part of the lemma it is enough to set
\begin{gather*}
f_{111}=b_{111} , \qquad f_{112}=b_{112}+\alpha S, \qquad f_{122}=b_{122}+\alpha T, \qquad f_{222}=b_{222}
\end{gather*}
and the corresponding $\E_{\D}$, with $\D$ given by~\eqref{eq.D.new}, is precisely described by~\eqref{eqGoursat} with $f_{112}=f_{211}$ and $f_{122}=f_{212}$.
\end{proof}

\begin{Corollary}\label{CorollarioCheSeFunzionaSarebbeCarino}
 Let $\E$ be an equation with $\MM^\E=\P\D_1$. Then $\E=\E_{\D_1}$.
\end{Corollary}

\begin{proof}
Observe that any $L_{m^2}$, with $m^2\in\E$, always contains a characteristic line $\elle$ correspon\-ding to a f\/ixed linear factor of the symbol, i.e., an element $\elle$ belonging to $\MM^\E$ (see Section \ref{subsubComponenteLineare}). In other words, inclusion \eqref{eqcuoricino} can be made more precise:
\begin{gather}\label{equazioneSempreVeraPerFortuna}
\E_{\MM^\E}\supseteq \E.
\end{gather}
 It remains to prove that the inverse of inclusion \eqref{equazioneSempreVeraPerFortuna} is valid when $\MM^\E=\P\D_1$.
 Indeed, Lemma~\ref{LemEd1QLMonge} and Lemma~\ref{Lemma-ex-Considerazioni-Forma-Generale} together guarantee that the left-hand side of~\eqref{equazioneSempreVeraPerFortuna} is a closed submanifold (possibly with singularities, see Section~\ref{sezChiarificatrice}) of codimension one.
But the right-hand side of~\eqref{equazioneSempreVeraPerFortuna} is a closed submanifold (again with possible singularities) of the same dimension of the left-hand side. Hence, the two of them must coincide as well.
\end{proof}

\subsubsection{Proof of the statement~1 of Theorem \ref{thFormMag}}\label{SubStat1}
The proof of its f\/irst claim is contained in the next corollary.
\begin{Corollary}\label{corollariocorollarioso}
A $3\Rd$ order MAE $\E_\D$ is quasi-linear if and only $\dim\D^v=2$.
\end{Corollary}
\begin{proof}
 A direct consequence of Lemmas \ref{LemQLsonoMonge}, \ref{LemEd1QLMonge} and \ref{Lemma-ex-Considerazioni-Forma-Generale}.
\end{proof}

\looseness=-1
The remainder of statement 1 is concerned with the ``other component'' $\MMM^\E$ of $\M^\E$, i.e., the one associated with the quadratic factor of $\Ann V\E$ (see Section \ref{subsubComponenteLineare}). Recall that $\MMM^\E$ might be empty, in which case there is nothing to prove. At the far end, there is the case when $\MMM^\E$ is, in its turn, decomposable, i.e., it consists of two linear sheets, which is dealt with by Lemma \ref{LemmaCalcolosoByGiovanni} below.

\begin{Lemma}\label{LemmaCalcolosoByGiovanni}
 Let $\E=\{F=0\}$ be a $3\Rd$ order quasi-linear PDE such that $\operatorname{Smbl} F$ is completely decomposable, i.e.,
 \begin{gather}\label{eqSymbDec}
\operatorname{Smbl} F=(\xi+k_1\eta)\odot(\xi+k_2\eta)\odot(\xi+k_3\eta).
\end{gather}
If
 \begin{gather}\label{eqDi}
\D_i:=\Span{\left. D_1\right|_{m^2}+k_i\left. D_2\right|_{m^2}\,|\, m^2\in \E_{m^1}}
\end{gather}
is the $3D$ sub-distribution of $\C^1$ corresponding\footnote{See also, on this concern, the construction of $\MM^\E$ in Section~\ref{subsubComponenteLineare}.} to the $i\Th$ linear factor of~\eqref{eqSymbDec} and
\begin{gather}\label{eqDiPrimo}
\overline{\D}=\Span{\overline{h},\overline{\D}^v}\subseteq\C^1
\end{gather}
is a generic $3D$ sub-distribution of $\C^1$ with $\dim \overline{\D}^v=2$, then $\overline{\D} $
 satisfies condition~\eqref{eqFormulaMagica} if and only if $\overline{\D}=\D_j$ with $j\neq i$.
\end{Lemma}

\begin{proof}
To begin with, \eqref{eqSymbDec} dictates some restrictions on $F$, which must be of the form 
 \begin{gather}\label{eqFdec}
F=p_{111}+(k_1+k_2+k_3)p_{112}+(k_1k_2+k_1k_3+k_2k_3)p_{122}+k_1k_2k_3p_{222}+c.
\end{gather}
We begin with the ``homogeneous'' case, i.e., we assume $c=0$ since, as we shall see at the end of the proof, the general case can be easily brought back to this one.

Equating \eqref{eqFdec} to zero allows to express $p_{111}$ as a linear combination of $p_{112}$, $p_{122}$, and $p_{222}$, i.e., to identify $\E_{m^1}$ with $\R^3\equiv\{(p_{112},p_{122},p_{222})\}$. In turn, this makes it possible to parametrize the space of vertical elements of \eqref{eqDi} by three real parameters, viz.,
\begin{gather}
 \D_i^v= \big\langle{-}({ (k_1+k_2+k_3-k_i)p_{112}+(k_1k_2+k_1k_3+k_2k_3)p_{122}+ k_1k_2k_3p_{222}})\partial_{p_{11}} \nonumber\\
\hphantom{\D_i^v= \big\langle}{} +({ p_{112}+k_ip_{122}})\partial_{p_{12}}+({ p_{122}+k_ip_{222}})\partial_{p_{22}}\,|\, { p_{112},p_{122},p_{222}}\in\R\big\rangle. \label{eqDiPar}
\end{gather}
It is worth observing that, in compliance with Lemma \ref{LemQLsonoMonge}, the dimension of $\D_i^v$, which
equals the rank of the $3\times 3$ matrix
\begin{gather}\label{eqMi}
M_i=\left|\begin{matrix}-k_1-k_2-k_3+k_i & -k_1k_2-k_1k_3-k_2k_3 & -k_1k_2k_3 \\1 & k_i & 0 \\0 & 1 & k_i\end{matrix}\right|
\end{gather}
 is 2: indeed,
$\det M_i= -(k_1-k_i)(k_2-k_i)(k_3-k_i)$
vanishes for $i=1,2,3$.

In order to f\/ind a basis for~\eqref{eqDiPar}, regard the matrix~\eqref{eqMi} as a (rank-two) homomorphism $M_i\colon \R^3\longrightarrow V=\Span{\partial_{p_{11}},\partial_{p_{12}},\partial_{p_{22}}}$ and compute its kernel:
\begin{gather*}
\ker M_i=\Span{k_i^2,-k_i,1}.
\end{gather*}
Then, independently on~$i$ (and on the value of~$k_i$ as well),
$\R^3=\Span{(1,0,0), (0,1,0)}\oplus \ker M_i$,
and
$\D_i^v=\langle M_i\cdot (1,0,0), M_i\cdot (0,1,0)\rangle$
 is the sought-for basis. In other words, instead of~\eqref{eqDi}, we shall work with the handier description
$\D_i=\Span{h_i, \D_i^v}$,
where
\begin{gather}
h_i = D_1 + k_iD_2 ,\nonumber\\
\D_i^v = \Span{ -(k_1+k_2+k_3-k_i) \partial_{p_{11}} + \partial_{p_{12}}, -(k_1k_2+k_1k_3+k_2k_3)\partial_{p_{11}}+k_i \partial_{p_{12}}+ \partial_{p_{22}} }.\label{eqDiV}
\end{gather}
Concerning \eqref{eqDiPrimo}, introduce similar descriptions of its horizontal and vertical part:
\begin{gather}
\overline{h} = aD_1+bD_2,\qquad
\overline{\D}^v = \Span{r_1\partial_{p_{11}}+s_1 \partial_{p_{12}}+t_1 \partial_{p_{22}},r_2\partial_{p_{11}}+s_2 \partial_{p_{12}}+t_2 \partial_{p_{22}}}.\label{eqVprimo}
\end{gather}
This concludes the preliminary part of the proof. Now impose condition \eqref{eqFormulaMagica}:
\begin{gather}\label{eqRicetteMagiche}
\Omega\big(h_i,\overline{\D}^v\big)=\Omega\big(\overline{h},\D_i^v\big)=\textrm{1D subspace}.
\end{gather}
Observe that \eqref{eqRicetteMagiche} consists, in fact, of two requirements, which are going to be dealt with separately.

The f\/irst one corresponds to the one-dimensionality of the subspace
\begin{align}
 \Omega\big(h_i,\overline{\D}^v\big)&=\Span{r_1\partial_{p_{1}}+s_1 \partial_{p_{2}}+k_i(s_1 \partial_{p_{1}}+t_1 \partial_{p_{2}}),r_2\partial_{p_{1}}+s_2 \partial_{p_{2}}+k_i(s_2 \partial_{p_{1}}+t_2 \partial_{p_{2}})}\nonumber\\
 &=\Span{(r_1+k_is_1)\partial_{p_{1}}+(s_1+k_it_1) \partial_{p_{2}}, (r_2+k_is_2)\partial_{p_{1}}+(s_2+k_it_2) \partial_{p_{2}}},\label{eqOmegaHDiPrimo}
\end{align}
i.e., to the equation
\begin{gather}\label{eqMostruosaChePerFortunaNonServe}
\det \left|\begin{matrix}r_1+k_is_1 & s_1+k_it_1 \\r_2+k_is_2 & s_2+k_it_2\end{matrix}\right| = k_i^2 A -k_i B + C =0,
\end{gather}
with
\begin{gather}\label{eqCoeffA-C}
A= \det\left|\begin{matrix}s_1 & t_1 \\s_2 & t_2\end{matrix}\right|, \qquad 
B= - \det\left|\begin{matrix}r_1 & t_1 \\r_2 & t_2\end{matrix}\right|,\qquad 
C= \det\left|\begin{matrix}r_1 & s_1 \\r_2 & s_2\end{matrix}\right|. 
\end{gather}
Interestingly \looseness=-1 enough, above coef\/f\/icients \eqref{eqCoeffA-C} characterize the dual direction of $\overline{\D}^v$, i.e., a non-zero co-vector $\omega'=Adp_{11}+Bdp_{12}+Cdp_{22}\in V^\ast$, def\/ined up to a non-zero constant, such that
\begin{gather}\label{Dvprimodual}
\overline{\D}^v=\ker \omega'.
\end{gather}
This dual perspective on $\overline{\D}^v$ allows to rewrite \eqref{eqMostruosaChePerFortunaNonServe} as
\begin{gather}\label{eqMostruosaChePerFortunaNonServeMigliore}
\omega'\big(k_i^2\partial_{p_{11}}-k_i\partial_{p_{12}}+ \partial_{p_{22}}\big)=0.
\end{gather}
In turn, \eqref{eqMostruosaChePerFortunaNonServeMigliore} dictates the form of
 $\omega'$:
\begin{gather}\label{eqFormaDellaFOrma}
\omega'=Adp_{11}+Bdp_{12}+\big(Bk_i-Ak_i^2\big)dp_{22}, \qquad (A,B)\neq (0,0).
\end{gather}
Together, \eqref{eqFormaDellaFOrma} and \eqref{Dvprimodual} produce a simplif\/ied expression
\begin{gather}\label{eqVprimoPiuBella}
 \overline{\D}^v=\Span{\big(Ak_i^2-Bk_i\big)\partial_{p_{11}}+A\partial_{p_{22}}, -B\partial_{p_{11}}+A\partial_{p_{12}}},
\end{gather}
which, in comparison with \eqref{eqVprimo}, depending on 6 parameters, needs only~2 of them, or even~1, if a non-zero constant is neglected.
Now, thanks to~\eqref{eqVprimoPiuBella}, it is easier to see that the space~\eqref{eqOmegaHDiPrimo} is~1D. Indeed,~\eqref{eqOmegaHDiPrimo} reads
\begin{gather*}
 \Omega\big(h_i,\overline{\D}^v\big) = \Span{\big(Ak_i^2-Bk_i\big)\partial_{p_{1}}+Ak_i\partial_{p_{2}},-B\partial_{p_{1}}+A\partial_{p_{2}}+Ak_i\partial_{p_{1}}},
\end{gather*}
where the f\/irst vector equals the second multiplied by $k_i$, so that it can be further simplif\/ied:
\begin{gather}\label{eqCheServiraDopo}
 \Omega(h_i, {\D}^v) = \Span{( Ak_i-B)\partial_{p_{1}}+A\partial_{p_{2}} }.
\end{gather}
We can pass to the other condition dictated by \eqref{eqRicetteMagiche}. In particular, the subspace
\begin{gather}
 \Omega\big(\overline{h},\overline{\D}^v\big) = \big\langle a (-(k_1+k_2+k_3-k_i) \partial_{p_{1}} + \partial_{p_{2}}) + b \partial_{p_{1}},\nonumber\\
\hphantom{\Omega\big(\overline{h},\overline{\D}^v\big) = \big\langle}{}
 (-(k_1k_2+k_1k_3+k_2k_3)\partial_{p_{1}}+k_i \partial_{p_{2}}) + b ( k_i \partial_{p_{1}}+ \partial_{p_{2}})\big\rangle\nonumber\\
\hphantom{\Omega\big(\overline{h},\overline{\D}^v\big) }{} = \big\langle (b-a(k_1+k_2+k_3-k_i)) \partial_{p_{1}} + a\partial_{p_{2}} ,\nonumber\\
\hphantom{\Omega\big(\overline{h},\overline{\D}^v\big) = \big\langle}{} (bk_i-a(k_1k_2+k_1k_3+k_2k_3))\partial_{p_{1}}+(ak_i+b) \partial_{p_{2}} \big\rangle\label{eqPrimoGeneratore}
\end{gather}
must be 1D, i.e.,
\begin{gather}\label{eqConProduttoria}
\det \left|\begin{matrix}b-a(k_1+k_2+k_3-k_i) & a \\bk_i-a(k_1k_2+k_1k_3+k_2k_3) & ak_i+b\end{matrix}\right|=\prod_{j\neq i} (b-ak_j)=0.
\end{gather}
Above equation \eqref{eqConProduttoria}, makes it evident that, for any $i=1,2,3$, the space~\eqref{eqCheServiraDopo} is 1D if and only if
\begin{gather}\label{eqSempliceCheNonMiAspettavo}
b=ak_j,\qquad j\in\{1,2,3\}{\setminus}\{i\}, \qquad a\neq 0,
\end{gather}
which corresponds to
\begin{gather}\label{eqAccPrimo}
\overline{h}=D_1 +k_jD_2 ,\qquad j\in\{1,2,3\}{\setminus}\{i\}.
\end{gather}
Plugging \eqref{eqSempliceCheNonMiAspettavo} into \eqref{eqPrimoGeneratore}, we f\/ind a unique vector
\begin{gather}\label{unicoGenComm2}
(k_i+k_j-k_1-k_2-k_3) \partial_{p_{1}} + \partial_{p_{2}}
\end{gather}
generating $ \Omega(\overline{h}, {\D}^v) $. Introducing the complement $c(i,j)$ of $\{i,j\}$ in $\{1,2,3\}$, \eqref{unicoGenComm2} reads
\begin{gather*}
 \Omega\big(\overline{h}, {\D}^v\big) = \Span{-k_{c(i,j)}\partial_{p_{1}} + \partial_{p_{2}} },\qquad j\in\{1,2,3\}{\setminus}\{i\}.
\end{gather*}
To conclude the proof, recall that, besides their one-dimensionality, \eqref{eqRicetteMagiche} also requires the equa\-li\-ty of the subspaces \eqref{eqVprimoPiuBella} and \eqref{eqCheServiraDopo}, i.e.,
$ \Omega(h_i,\overline{\D}^v)=\Span{( Ak_i-B)\partial_{p_{1}}+A\partial_{p_{2}} , -k_{c(i,j)}\partial_{p_{1}} + \partial_{p_{2}} }$ $= \Omega(\overline{h}, {\D}^v)$
or, alternatively,
 \begin{gather}\label{codizioneUguaglainzaCommutaori}
\det\left|\begin{matrix} Ak_i-B & A \\ -k_{c(i,j)} & 1 \end{matrix}\right| =0.
\end{gather}
Thanks to \eqref{codizioneUguaglainzaCommutaori},
 \begin{gather}\label{eqFINALMENTE}
B=A(k_i+k_{c(i,j)}),\qquad j\in\{1,2,3\}{\setminus}\{i\},
\end{gather}
and \eqref{eqFINALMENTE} allows to eliminate $B$ from \eqref{eqVprimoPiuBella}:
\begin{gather}\label{eqVprimoPiuBellaAncora}
\overline{\D}^v= \Span{\big(Ak_i^2-A(k_i+k_{c(i,j)})k_i\big)\partial_{p_{11}}+A\partial_{p_{22}}, -A(k_i+k_{c(i,j)})\partial_{p_{11}}+A\partial_{p_{12}}}.
\end{gather}
Being a non-zero constant\footnote{If $A=0$, then from \eqref{eqVprimoPiuBella} follows $\dim\overline{D}^v=1$.}, $A$ can also be removed from \eqref{eqVprimoPiuBellaAncora}, which becomes:
 \begin{gather}\label{eqVprimoPiuBellaAncoraDavvero}
\overline{\D}^v= \Span{ -k_{c(i,j)}k_i\partial_{p_{11}}+\partial_{p_{22}}, (k_{c(i,j)}+k_i)\partial_{p_{11}}+\partial_{p_{12}}}.
\end{gather}
To enlighten the conclusions, it is useful to rewrite together \eqref{eqAccPrimo} and \eqref{eqVprimoPiuBellaAncoraDavvero} above: there are exactly~2 distributions~$\overline{\D}$ which are ``compatible'' (in the sense of~\eqref{eqRicetteMagiche}) with the distribution~$\D_i$ given by~\eqref{eqDi}. More precisely, their horizontal and vertical parts are
\begin{gather}
\overline{h} = D_1 +k_jD_2,\qquad
\overline{\D}^v = \Span{ -k_{c(i,j)}k_i\partial_{p_{11}}+\partial_{p_{22}}, (k_{c(i,j)}+k_i)\partial_{p_{11}}-\partial_{p_{12}}},\label{eqDiPrimoRipetuta}
\end{gather}
respectively, for the only two possible values of $ j\in \{1,2,3\}{\setminus}\{i\}$. One only needs to reali\-ze that~\eqref{eqDiPrimoRipetuta} is one of the $\D_i^v$'s from~\eqref{eqDiV}.
To this end, rewrite~\eqref{eqDiV} replacing $i$ with $j$:
$\D_j^v=\langle -(k_1+k_2+k_3-k_j) \partial_{p_{11}} + \partial_{p_{12}}, -(k_1k_2+k_1k_3+k_2k_3)\partial_{p_{11}}+$ $k_j \partial_{p_{12}}- \partial_{p_{22}} \rangle$,
and subtract from the second vector the f\/irst one multiplied by $k_j$:
\begin{gather}
\D_j^v=\big\langle {-}(k_1+k_2+k_3-k_j) \partial_{p_{11}} \nonumber\\
\hphantom{\D_j^v=\big\langle}{} + \partial_{p_{12}}, (k_j(k_1+k_2+k_3-k_j)-(k_1k_2+k_1k_3+k_2k_3))\partial_{p_{11}} + \partial_{p_{22}} \big\rangle.\label{eqNonNEPosooPiu}
\end{gather}
Formulas
\begin{gather*}
k_j(k_1+k_2+k_3-k_j) = k_j\sum_{l\neq j}k_l, \qquad
k_1k_2+k_1k_3+k_2k_3 = k_j\sum_{l\neq j}k_l +\prod_{l\neq j}k_l,
\end{gather*}
show that the coef\/f\/icient of $\partial_{p_{11}} $ in the second vector of \eqref{eqNonNEPosooPiu} reduces to $-\prod\limits_{l\neq j}k_l$. Hence, \eqref{eqNonNEPosooPiu} reads
\begin{gather}\label{eqDjvConSommatorieEProduttorie}
\D_j^v=\bigg\langle{-}\bigg(\prod_{l\neq j}k_l \bigg)\partial_{p_{11}} + \partial_{p_{22}},\bigg( \sum_{l\neq j}k_l\bigg) \partial_{p_{11}} - \partial_{p_{12}} \bigg\rangle.
\end{gather}
Comparing now \eqref{eqDiPrimoRipetuta} with \eqref{eqDjvConSommatorieEProduttorie} it is evident that $\overline{\D}^v$ must equal $\D_j^v$, with $j\neq i$, and the proof of the ``homogeneous'' case is complete.

To deal with the general case, denote by $\E$ the equation determined by $F$ as in \eqref{eqFdec}, with $c=0$, and by $\E_c$ a ``inhomogeneous'' equation, i.e., one with $c\neq 0$. Observe that there is a~natural identif\/ication $i_c\colon \Span{D_1,D_2}\longrightarrow\Span{D_1-c\partial_{p_{11}},D_2}$ between horizontal planes, giving rise to an automorphism $\varphi_c:=i_c\oplus\id_{VM^{(1)}}$ of $\C^1$. Easy computations show that $\varphi_c(\D_i)$, $i=1,2,3$, are precisely the three distributions associated with the factors of the symbol of $\E_c$, and plainly $\Omega(h,\D^v)=\Omega(\varphi_c(h),\varphi_c(\D^v))$. So, the ``inhomogeneous'' case reduces to the ``homogeneous'' one, which has been established above.
\end{proof}

\begin{Corollary}[proof of statement~1 of Theorem~\ref{thFormMag}]\label{corRicettaMagica}
 Let $\E=\{F=0\}$ be a quasi-linear $3\Rd$ order PDE, and $U\subseteq M^{(1)}$ the open
 locus where the bundle $\MMM ^\E$ is not empty. Then
$\M^\E =\bigcup_{i=1}^3\P\D_i$ on~$U$,
 where $\D_i\subseteq\C|_U$ are {orthogonal} $($see Definition {\rm \ref{defTriOrt})} with respect to the meta-symplectic structure on $M^{(1)}$ and $\E=\E_{\D_i}$ on $U$, for all $i=1,2,3$.
\end{Corollary}
If $\E=\{F=0\}$ is a quasi-linear $3\Rd$ order PDE, and $m^1\in M^{(1)}$ is such that $\left(\MMM^\E\right)_{m_1}\neq\varnothing$, then it means that the function $\Delta$ (see Remark \ref{RemTipiDiEquazioni}) is positive on the whole f\/ibre $M^{(1)}_{m^1}$ and, hence, on the neighbouring f\/ibres. So, $U$ is indeed open.
\begin{proof}[Proof of Corollary \ref{corRicettaMagica}.]
 Let $m^1\in U$. Then $\MMM ^\E\neq \varnothing$, i.e., there is a point $m^2\in \E_{m^1}$ such that $\operatorname{Smbl}_{m^2}F$ is fully decomposable. But, in view of the quasi-linearity of $F$, this means that~$\operatorname{Smbl} F$ is fully decomposable on the whole f\/ibre $\E_{m^1}$, i.e., $\E$ fullf\/ills the hypotheses of Lemma~\ref{LemmaCalcolosoByGiovanni} on~$U$.
\end{proof}

\subsubsection{Proof of the statement 2 of Theorem \ref{thFormMag}}\label{SubStat2}

The f\/irst part of this statement merely rephrases the f\/ist part of statement 1.

Proving its second part consists in completing the proof of Lemma~\ref{LemQLsonoMonge} above with the computation of $\MMM^\E$, which corresponds to the quadratic factor (see Section~\ref{subsubComponenteLineare}) of~\eqref{eqSimboloPazzesco}. Formally, the procedure to get a description of $\MMM^\E$ out of the quadratic factor of \eqref{eqSimboloPazzesco} is the same as that for obtaining \eqref{eq.V1.Goursat.nonlin} from its linear factor, but it is more involved when it comes to computations, and this is the reason why we kept them separate.

So, suppose that $\MMM^\E\neq\varnothing$. This means that there are points where the quadratic factor of~\eqref{eqSimboloPazzesco} admits real roots.
If $\Delta$ is the discriminant of such quadratic factor, then one its factor is given by
\begin{gather}\label{eqLineaNonLineare}
\xi_1-\frac{T(p_{112}-f_{211})-R(p_{222}-f_{222})+\sqrt{\Delta}}{2(T(p_{122}-f_{212})-S(p_{222}-f_{222}))}\xi_2,
\end{gather}
and the corresponding line $\elle_1$ is obtained by replacing\footnote{This is a circumstance where the $\xi_i$'s and the $D_i$'s are vector f\/ields, not vectors (still, this notation is consistent with the conventions established in Remark~\ref{remcoordELLE}).} $\xi_i$ with $D_i|_\E$ in \eqref{eqLineaNonLineare}.
Now it is possible to assign certain values (namely $(0,0,0)$, $(1,0,0)$, $(0,1,0)$, $(0,0,1)$) to the internal parameters of~$\E$ (which are $p_{112}$, $p_{122}$, $p_{222}$), in such a way that the corresponding lines $\elle_1$ do not lie all in the same 3D subspace of~$\C^1$.
This shows that~$\MMM^\E$ cannot contain any linear irreducible component.

In order to prove that
\begin{gather*}
\E_{\D_1}=\E_{\MMM^\E},
\end{gather*}
observe that the ``$\supseteq$'' inclusion follows from the fact that, by def\/inition, $\MMM^\E$ is made of characte\-ris\-tic lines of $\E_{\D_1}$ which are also strongly characteristic by Proposition \ref{proposizioneNotturna}, so that all Lagrangian planes containing a line in $\MMM^\E$ belong also to $\E_{\D_1}$ (see also Def\/inition \ref{defStrongChar}).

Inclusion ``$\subseteq$'' can be obtained by a straightforward computation: indeed the $3\times 5$ matrix whose f\/irst line is formed by the components of $\elle_1$ and the remaining ones by the components of the vectors spanning a generic Lagrangian plane belonging to $\E_{\D_1}$, is of rank $2$. The same result is attained by considering the line $\elle_2$ obtained by changing the sign of $\sqrt{\Delta}$ in \eqref{eqLineaNonLineare}.

\subsubsection{Proof of the statement 3 of Theorem \ref{thFormMag}}\label{sec.3.3}

Given any 3D sub-distribution $\D\subseteq \C^1$,
and the corresponding Goursat-type equation $\E_\D$ (see~\eqref{eqED}), we are now in position to clarify how many and which sort of elements can be found in the set
\begin{gather}\label{eqInsiemeDelleSorprese}
\big\{\widetilde{\D}\,|\, \widetilde{\D} \textrm{ is a 3D sub-distribution of }\C^1\textrm{ such that }\E_{\widetilde{\D}}=\E_\D\big\}.
\end{gather}
More precisely,
\begin{enumerate}\itemsep=0pt
\item[1)] in the points of $M^{(1)}$ where $\dim\D^v=3$, we have $\E_\D=M^{(2)}$, so that the latter cannot be even considered as a $3\Rd$ order PDE according to our understanding of PDEs (see Section~\ref{sec.definitions.contact.and.prol}), and this case must be excluded \emph{a priori}\footnote{Still, the description of the set \eqref{eqInsiemeDelleSorprese} is very easy: it contains only $\D$, which in turn coincides with $VM^{(1)}$.};
\item[2)] in the points of $M^{(1)}$ where $\dim\D^v=1$, the equation $\E_\D$ is non-linear and the set \eqref{eqInsiemeDelleSorprese} contains only $\D$;
\item[3)] in the points of $M^{(1)}$ where $\dim\D^v=2$, the equation $\E_\D$ is quasi-linear and the set \eqref{eqInsiemeDelleSorprese} contains three elements $\D_i$, $i=1,2,3$, possibly repeated and comprising $\D$ itself, which
are {orthogonal} each other (see Def\/inition~\ref{defTriOrt}).
\end{enumerate}

\begin{proof}[Proof of 2 and 3]
 Let $\widetilde{\D}$ be an element of the set \eqref{eqInsiemeDelleSorprese}. Then, in particular,
 \begin{gather}\label{eqConiUguali}
\M^{\E_\D}=\M^{\E_{\widetilde{\D}}}.
\end{gather}
In view of Corollary \ref{ColloarioCheServiraPrimaOPoi}, equality \eqref{eqConiUguali} implies
 \begin{gather}\label{eqInclyusoneIncorciata1-2}
\P\D \subseteq \M^{\E_{\widetilde{\D}}},\qquad 
\P\widetilde{\D} \subseteq \M^{\E_{\D}}. 
\end{gather}
In other words, condition \eqref{decCono} is fulf\/illed by both equations $\E_{\widetilde{\D}}$ and $ \E_\D$, so that the statements~1 and~2 of Theorem~\ref{thFormMag} (see above Sections~\ref{SubStat1} and~\ref{SubStat2}) can be made use of. To this end, rewrite~\eqref{eqInclyusoneIncorciata1-2} as
 \begin{gather}
 \M^{\E_{\widetilde{\D}}} = \P\D \cup \MMM^{\E_{\widetilde{\D}}},\label{eqInclyusoneIncorciata1modif}\\
 \M^{\E_{D}} = \P\widetilde{\D} \cup \MMM^{\E_{{\D}}},\label{eqInclyusoneIncorciata2modif}
\end{gather}
respectively.
Because of statement 1 (reps., 2) of Theorem \ref{thFormMag}, from \eqref{eqInclyusoneIncorciata1modif} it follows that $\E_{\widetilde{\D}}$ is quasi-linear (resp., non-linear) if and only if $\dim\D^v=2$ (resp., 1) and from \eqref{eqInclyusoneIncorciata2modif} it follows that $\E_{\D}$ is quasi-linear (resp., non-linear) if and only if $\dim{\widetilde{\D}}^v=2$ (resp., 1). Summing up, $\dim \D^v=\dim\widetilde{\D}^v=1,2 $, and the two cases can be treated separately.

Directly from \eqref{eqInclyusoneIncorciata1modif} and \eqref{eqInclyusoneIncorciata2modif} we get
 \begin{gather}\label{eqInclyusoneIncorciata1-2modiff}
\P\widetilde{\D} \cup \MMM^{\E_{\widetilde{\D}}} = \P\D \cup \MMM^{\E_{\widetilde{\D}}},\qquad 
\P {\D} \cup \MMM^{\E_{ {\D}}} = \P\widetilde{\D} \cup \MMM^{\E_{{\D}}}. 
\end{gather}
Prove now~2. If $\dim \D^v=1$, then also $\dim\widetilde{\D}^v=1$ and the statement~2 of Theorem~\ref{thFormMag} (see
Section~\ref{SubStat2}) guarantees that $\MMM^{\E_{\widetilde{\D}}}$ (resp., $\MMM^{\E_{{\D}}}$) cannot contain $\P\D$ (resp., $\P\widetilde{\D}$). It follows from~\eqref{eqInclyusoneIncorciata1-2modiff} that $\P\widetilde{\D}=\P\D$.

Finally prove~3. If $\dim \D^v=2$, then also $\dim\widetilde{\D}^v=2$ and, in view of the statement 1 of Theorem \ref{thFormMag} (see
Section~\ref{SubStat1}), $\MMM^{\E_{\widetilde{\D}}}$ (resp., $\MMM^{\E_{{\D}}}$) is either empty, in which case $\P\widetilde{\D}=\P\D$, or consists of two, possibly repeated, distributions ``orthogonal'' to $\widetilde{\D}$ (resp.,~$\D$).
\end{proof}

The result presented in this section mirrors the analogous result for classical multi-dimensional MAEs (see \cite[Theorem~1]{MR2985508}), but displays some new and unexpected features: the threefold multiplicity of the notion of orthogonality and the distinction of the cases according to the dimension of the vertical part.

\section[Intermediate integrals of Goursat-type $3\Rd$ order MAEs]{Intermediate integrals of Goursat-type $\boldsymbol{3\Rd}$ order MAEs}\label{secIntInt}

Theorem \ref{thFormMag} established a link between $3\Rd$ order MAEs of Goursat-type, i.e., non-linear PDEs of order three, and 3D sub-distributions $\D$ of $\C^1$, i.e., linear objects involving (at most) second-order partial derivatives. Besides its aesthetic value, such a perspective also allows to formulate concrete results concerning the existence of solutions, as the f\/irst integrals of $\D$ can be made use of in order to f\/ind intermediate integrals of $\E_\D$, along the same lines of the classical case (see \cite[Section~6.3]{MR2985508} and~\cite{Alekseevsky2014144}, where a more general concept of intermediate integral is exploited to explicitly construct solutions to~$2\Nd$ order parabolic MAEs). Proposition~\ref{propIntInt} is the main result of this last section, showing that, for a Goursat-type $3\Rd$ order MAE~$\E$, the notions of an intermediate integral of $\E$ and of a f\/irst integral of \emph{any} distribution $\D$ such that $\E_\D=\E$ are actually the same. To facilitate its proof, we deemed it convenient to introduce equations of the form $\E=\E_\M$ where $\M=\P\Ddue$, with, \emph{for the first time in this paper}, $\dim\Ddue=2$. By analogy with \eqref{eqED}, we still write $\E_\Ddue$ instead of $\E_{\P\Ddue}$, but we warn the reader that, unlike all the cases considered so far, $\E_\Ddue$ is actually a \emph{system} of two independent equations (see Lemma \ref{LemmaPartoritoConTantaSofferenza} below).

Recall that a function $f\in C^\infty(M^{(1)})$ determines, in the neighborhood of an its non-singular point, the hyperplane distribution $\ker d f$ on $M^{(1)}$, each of whose leaves is identif\/ied by a value $c\in\R$. Following the same notation as \cite{MR2985508}, we set
\begin{gather*}
 \big(M^{(1)}\big)_{f-c}:=\big\{m^1\in M^{(1)}\,|\, f\big(m^1\big)=c\big\} .
\end{gather*}

\begin{Definition}
A function $f\in C^\infty(M^{(1)})$ is an \emph{intermediate integral} of a PDE $\E\subseteq M^{(2)}$ if
 \begin{gather}\label{eqDefIntInt}
\big(M^{(1)}\big)_{f-c}^{(1)}\subseteq\E
\end{gather}
for all $c\in\R$.
\end{Definition}

Inclusion \eqref{eqDefIntInt} can be interpreted as follows: if $f=c$ def\/ines a $2\Nd$ order PDE (in the sense given in Section \ref{sec.definitions.contact.and.prol}), then any solution of such an equation is also a solution of the $3\Rd$ order equation~$\E$.

\begin{Lemma}\label{LemmaPartoritoConTantaSofferenza}
 Let $\Ddue\subset\C^1$ be a $2D$ non-vertical
 sub-distribution.
 Then $\E_\Ddue\longrightarrow M^{(1)}$ is a non-linear bundle of rank~$2$.
\end{Lemma}

\begin{proof}
To simplify the notations, let $\E:=\E_\Ddue$. Then,
f\/ix $m^1\in M^{(1)}$ and observe that
$
\E_{m^1}=\{m^2\in M^{(2)}\,|\, L_{m^2}\cap\Ddue_{m^1} \neq 0\}
$
is a 2D manifold. Indeed, $\E_{m^1}$ is\footnote{This step may require the restriction to an open and dense subset.} a rank-one bundle over the 1D manifold $\P\Ddue$, its f\/ibres being the prolongations of the lines lying in $\Ddue$ (see \cite[Proposition 2.5]{MR1670044}).

To prove non-linearity of $\E$, write down $\Ddue$ in local coordinates as
\begin{gather*}
\Ddue=\Span{a_iD_1+b_iD_2+R_i\partial_{p_{11}}+S_i\partial_{p_{12}}+T_i\partial_{p_{22}}\,|\, i=1,2}
\end{gather*} and observe that $\E$ is given by the vanishing of the f\/ive $4\times 4$ minors of the matrix
\begin{gather*}
\left|\begin{matrix}a_1 & b_1 & R_1 & S_1 & T_1 \\a_2 & b_2 & R_2 & S_2 & T_2 \\1 & 0 & p_{111} & p_{112} & p_{122} \\0 & 1 & p_{112} & p_{122} & p_{222}\end{matrix}\right|.
\end{gather*}
Tedious \looseness=1 computations show that the only cases when such minors are quasi-linear, is that \mbox{either} when $\Ddue$ fails to be two-dimensional or when $\Ddue$ is vertical, which are forbidden by the hypo\-theses.
\end{proof}

\begin{Corollary}\label{corollarioUnPoPesante}
 Let $\Ddue\subset\C^1$ be a $2D$ sub-distribution contained in $\ker d f$. Then for any \mbox{$m^1\in M^{(1)}$} there is a $m^2\in M^{(2)}_{m^1}$ such that $L_{m^2}\subset \ker d_{m^1}f$ and $L_{m^2}\cap \Ddue_{m^1}=0$.
\end{Corollary}

\begin{proof}
 Assume for contradiction that
 \begin{gather}\label{eqInclusione}
\big\{m^2\in M^{(2)}\,|\, L_{m^2}\subset \ker d_{m^1}f\big\}\subseteq \big\{m^2\in M^{(2)}\,|\, L_{m^2}\cap \Ddue_{m^1}\neq 0\big\}
\end{gather}
and set $c:=f(m^1)$. Then \eqref{eqInclusione} means that, over the point $m^1$,
 \begin{gather}\label{eqFinalmenteContraddittoria}
 \big(M^{(1)}\big)_{f-c}^{(1)}\subseteq \E_\Ddue ,
\end{gather}
which is impossible, since the left-hand side of~\eqref{eqFinalmenteContraddittoria} is linear while its right-hand side is either empty (in the case when $\Ddue\subseteq VM^{(1)}$) or non-linear, thanks to above Lemma~\ref{LemmaPartoritoConTantaSofferenza}.
\end{proof}

\begin{Proposition}\label{propIntInt}
 Let $\E$ be a Goursat-type $3\Rd$ order MAE and $f\in C^\infty(M^{(1)})$. Then
 the following statements are equivalent:
 \begin{itemize}\itemsep=0pt
 \item $f$ is an intermediate integral of $\E$;
\item there exists a $3D$ distribution $\D\subset\C^1$ such that $f$ is a first integral of $\D$ and $\E_\D=\E$.
\end{itemize}
\end{Proposition}

\begin{proof}
 If $f$ is a f\/irst integral of $\D$, then
 \begin{gather}\label{eqFLAG345}
\D\subseteq \C^1\cap \ker df \subseteq \C^1
\end{gather}
is a $(3,4,5)$-type f\/lag of distributions ($\C^1$ does not posses f\/irst integrals). By def\/inition,
\begin{gather}\label{eqDefProlungMf}
\big(M^{(1)}\big)_{f-c}^{(1)}=\big\{m^2\in M^{(2)}\,|\, f(\pi_{2,1}(m^2))=c,\ df|_{L_{m^2}}=0\big\}
\end{gather}
and condition $df|_{L_{m^2}}=0$ means precisely that the 2D subspace $L_{m^2}$ of $\C^1$ is also contained in $\ker df$, i.e., $L_{m^2}$ lies in the 4D subspace $\C^1\cap \ker df$ of~$\C^1$. Because of~\eqref{eqFLAG345}, the 3D subspace $\D$ is also contained in $\C^1\cap \ker df$ and, as such, it cannot fail to non-trivially intersect $L_{m^2}$. This means that~$m^2\in\E_\D$, and~\eqref{eqDefIntInt} follows from the arbitrariness of~$m^2$.

Suppose, conversely, that $f$ is an intermediate integral of $\E$. In view of \eqref{eqDefProlungMf} and inclusion~\eqref{eqDefIntInt}, one has that
\begin{gather}\label{eqImplicazioneInversaUnPoPezzottata}
L_{m^2}\subseteq\ker df\Rightarrow L_{m^2}\cap\D_{m^1}\neq 0
\end{gather}
for \emph{any} distribution $\D$ such that $\E_\D=\E$. We wish to show that, there is \emph{at least one} among these distributions such that $f$ is an its f\/irst integral, i.e., that the f\/irst inclusion of~\eqref{eqFLAG345} is satisf\/ied. Towards a contradiction, suppose that no distribution~$\D$ is contained in $\ker df$, i.e., that $\Ddue_f:=\D\cap\ker df$ is a 2D sub-distribution of the 4D distribution $\C^1\cap\ker df$, \emph{for all} distributions~$\D$ such that $\E_\D=\E$. Then, by Corollary~\ref{corollarioUnPoPesante}, it is possible to f\/ind a~(2D~subspace)~$L_{m^2}$ which is contained in~$\C^1\cap\ker df$ and also trivially intersects one~$\Ddue_f$, thus violating \eqref{eqImplicazioneInversaUnPoPezzottata}.
\end{proof}

\begin{Corollary}
 If the derived flag of $\D$ never reaches $TM^{(1)}$, then $\E_\D$ admits an intermediate integral.
\end{Corollary}

\subsection*{Acknowledgments}
The authors wish to express their gratitude towards the anonymous referees whose comments contributed to shape the paper into its f\/inal form.
The authors thank C. Ciliberto, E. Ferapontov and F. Russo for stimulating discussions. The research of the f\/irst author has been partially supported by the project ``Finanziamento giovani studiosi~-- Metriche proiettivamente equiva\-lenti, equazioni di Monge--Amp\`ere e sistemi integrabili'', University of Padova 2013--2015, by the project ``FIR (Futuro in Ricerca) 2013~-- Geometria delle equazioni dif\/ferenziali''. The research of the second author has been partially supported by the Marie Sk{\l}odowska--Curie Action No~654721 ``GEOGRAL'', by the University of Salerno, and by the project P201/12/G028 of the Czech Republic Grant Agency (GA \v{C}R). Both the authors are members of G.N.S.A.G.A.\ of I.N.d.A.M.

\pdfbookmark[1]{References}{ref}
\LastPageEnding


\begin{thebibliography}{99}
\footnotesize\itemsep=-0.8pt

\bibitem{MR1925763}
Agafonov S.I., Ferapontov E.V., Systems of conservation laws in the setting of
 the projective theory of congruences: reducible and linearly degenerate
 systems, \href{http://dx.doi.org/10.1016/S0926-2245(02)00105-5}{\textit{Differential Geom. Appl.}} \textbf{17} (2002), 153--173.

\bibitem{MR2985508}
Alekseevsky D.V., Alonso Blanco R., Manno G., Pugliese F., Contact geometry of
 multidimensional {M}onge--{A}mp\`ere equations: characteristics, intermediate
 integrals and solutions, \href{http://dx.doi.org/10.5802/aif.2686}{\textit{Ann. Inst. Fourier (Grenoble)}} \textbf{62}
 (2012), 497--524, \href{http://arxiv.org/abs/1003.5177}{arXiv:1003.5177}.

\bibitem{Alekseevsky2014144}
Alekseevsky D.V., Alonso Blanco R., Manno G., Pugliese F., Finding solutions of
 parabolic {M}onge--{A}mp\`ere equations by using the geometry of sections of
 the contact distribution, \href{http://dx.doi.org/10.1016/j.difgeo.2013.10.015}{\textit{Differential Geom. Appl.}} \textbf{33}
 (2014), suppl., 144--161.

\bibitem{MR2383541}
Alonso~Blanco R., Manno G., Pugliese F., Contact relative dif\/ferential
 invariants for non generic parabolic {M}onge--{A}mp\`ere equations,
 \href{http://dx.doi.org/10.1007/s10440-008-9204-8}{\textit{Acta Appl. Math.}} \textbf{101} (2008), 5--19.

\bibitem{MR2503974}
Alonso~Blanco R., Manno G., Pugliese F., Normal forms for {L}agrangian
 distributions on 5-dimensional contact manifolds, \href{http://dx.doi.org/10.1016/j.difgeo.2008.06.019}{\textit{Differential Geom.
 Appl.}} \textbf{27} (2009), 212--229, \href{http://arxiv.org/abs/0707.0683}{arXiv:0707.0683}.

\bibitem{ThesisMichi}
B{\"a}chtold M.J., Fold-type solution singularities and charachteristic
 varieties of non-linear PDEs, Ph.D.~Thesis, \href{http://dx.doi.org/10.5167/uzh-42719}{Universit\"at Z\"urich}, 2009.

\bibitem{MorBach}
B{\"a}chtold M.J., Moreno G., Remarks on non-maximal integral elements of the
 {C}artan plane in jet spaces, \href{http://dx.doi.org/10.1016/j.geomphys.2014.05.006}{\textit{J.~Geom. Phys.}} \textbf{85} (2014),
 185--195, \href{http://arxiv.org/abs/1208.5880}{arXiv:1208.5880}.

\bibitem{MR1670044}
Bocharov A.V., Chetverikov~V.N., Duzhin S.V., Khor'kova~N.G., Krasil'shchik~I.S., Samokhin~A.V., Torkhov~Yu.N., Verbovetsky~A.M., Vinogradov~A.M.,
 Symmetries and conservation laws for dif\/ferential equations of mathematical
 physics, \textit{Translations of Mathematical Monographs}, Vol.~182, Amer.
 Math. Soc., Providence, RI, 1999.

\bibitem{MR1139843}
Boillat G., Sur l'\'equation g\'en\'erale de {M}onge--{A}mp\`ere \`a plusieurs
 variables, \textit{C.~R.~Acad. Sci. Paris S\'er.~I Math.} \textbf{313}
 (1991), 805--808.

\bibitem{MR1194520}
Boillat G., Sur l'\'equation g\'en\'erale de {M}onge--{A}mp\`ere d'ordre
 sup\'erieur, \textit{C.~R.~Acad. Sci. Paris S\'er.~I Math.} \textbf{315}
 (1992), 1211--1214.

\bibitem{MR1083148}
Bryant R.L., Chern S.S., Gardner~R.B., Goldschmidt~H.L., Grif\/f\/iths~P.A.,
 Exterior dif\/ferential systems, \href{http://dx.doi.org/10.1007/978-1-4613-9714-4}{\textit{Mathematical Sciences Research
 Institute Publications}}, Vol.~18, Springer-Verlag, New York, 1991.

\bibitem{MR2805306}
De~Paris A., Vinogradov A.M., Scalar dif\/ferential invariants of symplectic
 {M}onge--{A}mp\`ere equations, \href{http://dx.doi.org/10.2478/s11533-011-0046-7}{\textit{Cent. Eur.~J. Math.}} \textbf{9}
 (2011), 731--751, \href{http://arxiv.org/abs/1102.0426}{arXiv:1102.0426}.

\bibitem{MR2369200}
De~Poi P., Mezzetti E., Congruences of lines in {${\mathbb P}^5$}, quadratic
 normality, and completely exceptional {M}onge--{A}mp\`ere equations,
 \href{http://dx.doi.org/10.1007/s10711-007-9228-7}{\textit{Geom. Dedicata}} \textbf{131} (2008), 213--230, \href{http://arxiv.org/abs/0710.5110}{arXiv:0710.5110}.

\bibitem{MR1397274}
Dubrovin B., Geometry of {$2$}{D} topological f\/ield theories, in Integrable
 Systems and Quantum Groups ({M}ontecatini {T}erme, 1993), \href{http://dx.doi.org/10.1007/BFb0094793}{\textit{Lecture
 Notes in Math.}}, Vol.~1620, Springer, Berlin, 1996, 120--348,
 \href{http://arxiv.org/abs/hep-th/9407018}{hep-th/9407018}.

\bibitem{MR1958112}
Ferapontov E.V., Decomposition of higher-order equations of {M}onge--{A}mp\`ere
 type, \href{http://dx.doi.org/10.1023/A:1022236709346}{\textit{Lett. Math. Phys.}} \textbf{62} (2002), 193--198,
 \href{http://arxiv.org/abs/nlin.SI/0205056}{nlin.SI/0205056}.

\bibitem{MR1504329}
Goursat E., Sur les \'equations du second ordre \`a {$n$} variables analogues
 \`a l'\'equation de {M}onge--{A}mp\`ere, \textit{Bull. Soc. Math. France}
 \textbf{27} (1899), 1--34.

\bibitem{MR1202431}
Kol{\'a}{\v{r}} I., Michor P.W., Slov{\'a}k J., Natural operations in
 dif\/ferential geometry, \href{http://dx.doi.org/10.1007/978-3-662-02950-3}{Springer-Verlag}, Berlin, 1993.

\bibitem{MR2813504}
Krasil'shchik J., Verbovetsky A., Geometry of jet spaces and integrable
 systems, \href{http://dx.doi.org/10.1016/j.geomphys.2010.10.012}{\textit{J.~Geom. Phys.}} \textbf{61} (2011), 1633--1674,
 \href{http://arxiv.org/abs/1002.0077}{arXiv:1002.0077}.

\bibitem{MR2389645}
Kruglikov B., Lychagin V., Geometry of dif\/ferential equations, in \href{http://dx.doi.org/10.1016/B978-044452833-9.50015-2}{Handbook of
 Global Analysis}, Elsevier Sci.~B.~V., Amsterdam, 2008, 725--771.

\bibitem{MR2352610}
Kushner A., Lychagin V., Rubtsov V., Contact geometry and non-linear
 dif\/ferential equations, \textit{Encyclopedia of Mathematics and its
 Applications}, Vol.~101, Cambridge University Press, Cambridge, 2007.

\bibitem{MR0067317}
Lax P.D., Milgram A.N., Parabolic equations, in Contributions to the Theory of
 Partial Dif\/ferential Equations, \textit{Annals of Mathematics Studies},
 Vol.~33, Princeton University Press, Princeton, N.~J., 1954, 167--190.

\bibitem{MR781584}
Lychagin V.V., Singularities of multivalued solutions of nonlinear dif\/ferential
 equations, and nonlinear phenomena, \href{http://dx.doi.org/10.1007/BF00580702}{\textit{Acta Appl. Math.}} \textbf{3}
 (1985), 135--173.

\bibitem{MR966202}
Lychagin V.V., Geometric theory of singularities of solutions of nonlinear
 dif\/ferential equations, \href{http://dx.doi.org/10.1007/BF01095431}{\textit{J.~Sov. Math.}} \textbf{51} (1990),
 2735--2757.

\bibitem{138544}
Moreno G., Submanifolds in the Grassmannian of $n$-dimensional subspaces
 determined by a~sub\-ma\-ni\-fold in the Grassmannian of $l$-dimensional subspaces,
 MathOverf\/low, 2013, available at \url{http://mathoverflow.net/q/138544}.

\bibitem{MorenoCauchy}
Moreno G., The geometry of the space of {C}auchy data of nonlinear {PDE}s,
 \href{http://dx.doi.org/10.2478/s11533-013-0292-y}{\textit{Cent. Eur.~J. Math.}} \textbf{11} (2013), 1960--1981,
 \href{http://arxiv.org/abs/1207.6290}{arXiv:1207.6290}.

\bibitem{MR1458653}
Morimoto T., Monge--{A}mp\`ere equations viewed from contact geometry, in
 Symplectic Singularities and Geometry of Gauge Fields ({W}arsaw, 1995),
 \textit{Banach Center Publ.}, Vol.~39, Polish Acad. Sci., Warsaw, 1997,
 105--121.

\bibitem{MR2279499}
Piccione P., Tausk D.V., The single-leaf {F}robenius theorem with applications,
 \textit{Resenhas} \textbf{6} (2005), 337--381, \href{http://arxiv.org/abs/math.DG/0510555}{math.DG/0510555}.

\bibitem{MR1372461}
Strachan I.A.B., On the integrability of a third-order {M}onge--{A}mp\`ere type
 equation, \href{http://dx.doi.org/10.1016/0375-9601(95)00944-2}{\textit{Phys. Lett.~A}} \textbf{210} (1996), 267--272.

\bibitem{LucaChar}
Vitagliano L., Characteristics, bicharacteristics and geometric singularities
 of solutions of {PDE}s, \href{http://dx.doi.org/10.1142/S0219887814600391}{\textit{Int.~J. Geom. Methods Mod. Phys.}} \textbf{11}
 (2014), 1460039, 35~pages, \href{http://arxiv.org/abs/1311.3477}{arXiv:1311.3477}.

\end{thebibliography}
\end{document}